\documentclass{article}
\usepackage{arxiv}
\pdfoutput=1
\usepackage[nocomma]{optidef}
\usepackage{array}
\usepackage{graphicx}
\usepackage{xcolor}
\usepackage{amsfonts}
\usepackage{relsize}
\usepackage{algorithm} 
\usepackage{algpseudocode} 
\usepackage{amsmath}
\usepackage{scalerel,tabstackengine,xpatch}
\setstacktabbedgap{1em}
\xpatchcmd\Centerstack{\strutlongstacks{T}}{}{}{}

\usepackage{tikz}
\usetikzlibrary{positioning}
\usetikzlibrary{shapes.geometric, arrows}

\newcommand{\exclude}[1]{}

\usepackage{natbib}

\begin{document}
\title{A Stochastic Lookahead Approach for Hurricane Relief Logistics Operations Planning under Uncertainty}

\author{
    Yanbin Chang \\
    Department of Industrial Engineering\\
    Clemson University\\ 
    Clemson, SC, USA\\
    \texttt{yanbinc@clemson.edu} \\
    \And
    Yongjia Song \\
    Department of Industrial Engineering\\
    Clemson University\\ 
    Clemson, SC, USA\\
    \texttt{yongjis@clemson.edu} \\
    \And
    Burak Eksioglu \\
    Department of Industrial Engineering\\
    University of Arkansas\\
    Fayetteville, AR, USA\\
    \texttt{burak@uark.edu} \\
}

\maketitle

\begin{abstract}
In the aftermath of a hurricane, humanitarian logistics plays a critical role in delivering relief items to the affected areas in a timely fashion. This paper proposes a novel stochastic lookahead framework that implements a two-stage stochastic programming model in a rolling horizon fashion to address the evolving uncertain logistics system state during the post-hurricane humanitarian logistics operations. The two-stage stochastic programming model that executes in this rolling horizon framework is formulated as a mixed-integer programming problem. The model aims to minimize the sum of transportation and social costs. The social cost is measured as a function of deprivation for unsatisfied demand. Our extensive numerical experiment results and sensitivity analysis demonstrate the effectiveness of the proposed approach in reducing the total cost incurred during the post-hurricane relief logistics operations compared to the two-stage stochastic programming model implemented in a static approach. 

\keywords{Stochastic programming \and Rolling horizon \and Disaster relief logistics \and Social cost}
\end{abstract}

\section{Introduction}
\label{intro}
Many regions in the United States and around the world are vulnerable to large-scale natural disasters such as earthquakes, hurricanes, and floods. Hurricanes, in particular, have been studied extensively because they affect large areas and populations and are relatively easier to predict compared to other disasters. Despite this fact, hurricane relief operations planning remains a challenging task. For instance, in 2005, Hurricane Katrina destroyed more than 200,000 houses and cost over \$100 billion \citep{deryugina_2018} although the National Hurricane Center provided accurate predictions of the strength and landfall location 48 hours before the storm hit \citep{knabb_2005}. In addition, even though an unprecedented number of critical resources was prepositioned based on the predictions, the Federal Emergency Management Agency (FEMA) was heavily criticized for failing to prepare enough relief items \citep{menzel_2006}. \citet{holguin_2012} pointed out that one of the most heavily criticized aspects was the inefficiency of the logistics operations, which did not deliver the relief items to the affected populations in a timely fashion after the hurricane.

The experience of Hurricane Katrina elucidates that a reliable and effective humanitarian relief logistics system is critical to reducing the negative impacts. \citet{van_2006} defines humanitarian relief logistics as the emergency operations involving mobilizing various resources that use skills and knowledge to help affected people in a disaster. The most challenging aspect in humanitarian relief operations is the coordination of logistics management in complicated and uncertain post-disaster environments \citep{richey_2009}. Compared to no-notice events such as earthquakes, it is more practicable to plan and manage the humanitarian relief logistics for advance-notice disasters such as hurricanes. Most emergency managers prefer a flexible and realistic optimization model to tackle the post-disaster humanitarian relief logistics problem \citep{vanajakumari_2016}. The main objective of this paper is to formulate and evaluate optimization models for the humanitarian relief logistics operations planning problem under various uncertain factors, such as the supply and demand for relief items. 

This paper addresses several gaps in humanitarian relief logistics operations. First, we provide a new method to calculate social costs (in terms of deprivation costs) incurred by the affected population under evolving demand and supply limit uncertainty. Second, we develop a two-stage stochastic programming model to deal with the uncertainty in the demand for and supply limit of relief items. Third, we propose a stochastic lookahead framework in the form of a rolling horizon approach to address hurricane relief logistics operations under evolving conditions. We show that the rolling horizon approach increases both the flexibility and the efficiency of hurricane relief logistics planning. 

The rest of the paper is organized as follows. We review the relevant literature in Section \ref{sec:2}. Section \ref{sec:3} describes the problem setting and our mathematical optimization formulations in detail. Section \ref{sec:4} presents our numerical results and sensitivity analysis. Finally, we conclude the paper and provide some future research directions in Section \ref{sec:5}. 

\section{Literature Review}
\label{sec:2}
Our study is related to four streams of literature: last mile delivery, social cost of disasters, stochastic programming, and rolling horizon. We provide brief reviews of the studies in these areas and explain how we contribute to each. 

\subsection{Last Mile Delivery} \label{sec:2.1}
In general, the last mile delivery problem focuses on scheduling efficient and effective deliveries of goods from local distribution centers to customers \citep{wang_2016}. The post-disaster last mile delivery of relief commodities to affected populations is a critical component of humanitarian relief logistics planning. \citet{haghani_1996} studied detailed routing plans for a humanitarian relief logistics problem with multiple transportation modes and various commodities. \citet{barbarosouglu_2002} developed a mathematical model for tactical and operational air transportation mission planning during a disaster. \citet{balcik_2008} proposed a two-phase modeling approach to address the delivery routes in a last-mile distribution system in a dynamic framework. In the first phase, they generated some candidate routes. Then they selected the optimal route for each period and determined the delivery amounts in the second phase. Our proposed model contains the following characteristics of the last mile distribution problem: (i) demand for relief items includes multiple commodities, (ii) both ground and air transportation are available for delivering the relief items, (iii) decision makers need to decide when and which potential temporary warehouse to open to facilitate the disaster relief logistics operations. Unlike many papers in the literature that include the vehicle routing aspect of logistics operations planning in their last mile delivery problem, we choose not to include it in our model because our focus is to tackle the stochastic nature of the problem. Although integrating the vehicle routing aspect of logistics operations in the overall planning model would make the model more realistic, it becomes much more computationally demanding, especially when we incorporate evolving uncertainty of system states over time. For this reason, we will focus on developing network flow type models for the last mile delivery problem in this paper. 

\subsection{Social Cost of Disasters}
\label{sec:2.2}
Social cost, primarily due to unmet demand, has been extensively modeled in the literature and incorporated in humanitarian relief logistics decision making. For example, \citet{perez_2016} pointed out that the objective functions in many humanitarian relief logistics studies minimize either (i) the logistics costs, (ii) a measure of human suffering, or (iii) a combination of both logistics costs and human suffering. Human suffering has been measured by the amount of unmet demand~\citep{ozdamar_2004, tzeng_2007, balcik_2008} or a  penalty-based function of unmet demand \citep{barbarosoglu_2004, chang_2007, yi_2007a, yi_2007, rawls_2010}. According to \citet{holguin_2013}, a limitation of using just the amount of unmet demand is that such models cannot distinguish the urgency level of certain demand. For example, areas that have been deprived of some small amount of a relief item for a long period of time should have higher priority over areas that have been deprived of a large amount of the same relief item for a much shorter period of time. Another limitation of models that measure human suffering by just the amount of unmet demand is that they assume demand backlog can be accumulated. However, this is not necessarily true for many relief commodities. For example, one cannot simply consume three days' worth of food and water when supplies arrive after being deprived for three days \citep{holguin_2013}.

Penalty-based models, on the other hand, can incorporate the urgency of unsatisfied demand as a function of the deprivation time. Yet, \citet{holguin_2013} suggested that penalty-based models are not ideal since they are unable to balance the operational issues and social welfare on a common scale. Thus, the subjective parameters can lead to biased solutions based on the preference structure of the researchers conducting the study. To overcome this problem, \citet{holguin_2013} formulated a model in which the objective is to minimize a total cost that consists of the logistics costs and the deprivation costs. To estimate the deprivation costs, they developed a function which maps the deprivation time (the amount of time that the demand goes unsatisfied) to the deprivation cost using certain socioeconomic characteristics. \citet{perez_2016} proposed a generic version of the deprivation cost function under the assumption that the demand remains constant over time. In this paper, we propose a new function to compute the deprivation cost that takes into account both the fluctuating demand and the deprivation time. 

\subsection{Two-Stage Stochastic Programming Models}
\label{sec:2.3}

Two-stage stochastic programming has been widely used to deal with uncertainty in the humanitarian relief logistics problem. In most of these models, the first-stage variables correspond to decisions regarding the inventory levels and the locations of facilities for stockpiling relief commodities, and the second-stage variables correspond to the logistics operational decisions that transport the relief commodities to the affected area after the disaster \citep{barbarosoglu_2004, mete_2010, van_Hentenryck_2010,  rawls_2010, hu_2015, tofighi_2016}. For example, \citet{salmeron_2010} proposed a two-stage stochastic programming model for evacuating critical population. In the first stage, decision variables correspond to the capacity extension of warehouses, medical facilities, and air transportation ramp spaces. The second stage decisions correspond to the logistics operations after the hurricane, which include the allocation of the relief commodities and deploying vehicles to rescue the critical population. \citet{li_2012} built a two-stage stochastic programming model that identifies the shelters to be maintained overtime in the first stage. In the second stage, they select shelters outside the affected area to open and prepare for evacuees. Table \ref{table:character_two-stage} summarizes the common characteristics of two-stage stochastic programming models in the literature used for disaster relief logistics problems.

\begin{table}[htbp]\caption{Summary of Two-stage Stochastic Programming Models Used in Disaster Relief Logistics Planning Problems.}\label{table:character_two-stage}
\setlength{\tabcolsep}{1pt}
\footnotesize
\begin{tabular}{>{\centering}p{0.28\textwidth}>{\centering}p{0.25\textwidth}>{\centering}p{0.25\textwidth}>{\centering\arraybackslash}p{0.2\textwidth}}
\hline
{} &1st-stage decisions &2nd-stage decisions &Uncertainty\\
\noalign{\smallskip}\hline
\citet{barbarosoglu_2004} &preposition inventory & logistics schedules &demand \& supply \\
\citet{mete_2010} &preposition location \& inventory & logistics schedules & demand \\
\citet{salmeron_2010} &expansion of facilities &logistics schedules &demand \\
\citet{van_Hentenryck_2010} &preposition location \& inventory &logistics schedules &supply \\
\citet{rawls_2010} &preposition location \& inventory &logistics schedules &supply \\
\citet{li_2012} &permanent shelter location & temporary shelter location &demand \\
\citet{hu_2015} &preposition location \& inventory &logistics schedules &demand \\
\citet{tofighi_2016} &preposition inventory &logistics schedules &demand \& supply \\
\hline
\end{tabular}
\end{table}

In many cases, a large number of scenarios are necessary to properly demonstrate uncertainty, which may lead to computational challenges for such two-stage stochastic programming models. Furthermore, introducing integer variables (\textit{e.g.}, opening of warehouses) into the stochastic programs potentially increases the complexity of the models. Although many decomposition methods (\textit{e.g.}, the integer L-shaped approach or branch-and-cut algorithms) are typically used to improve the computing efficiency, such models with a large number of scenarios are still computationally challenging. One of our key contributions is that the proposed stochastic lookahead framework, which embeds two-stage stochastic programming models in a rolling-horizon approach, could produce high-quality decisions without resorting to a large number of scenarios, according to our extensive computational results. This makes the implementation of the proposed framework relatively more efficient computationally than alternative offline approaches. 

\subsection{The Rolling Horizon Approach}
\label{sec:2.4}
The rolling horizon approach is a common optimization based approach for solving multi-period planning problems that captures the changing values of problem parameters~\citep{hasani_2018}. \citet{rivera_2016} suggested that it is possible to run an optimization model using a rolling horizon approach, whose parameters need to be continuously updated as more system information becomes available. They built a model to serve demand according to the level of urgency of demand points given by the continuously updating demand information during disaster relief operations. \citet{vanajakumari_2016} presented an integrated humanitarian relief logistics model in a multi-period setting that determines the optimal locations and inventory levels of prepositioned supplies and the routing plans of trucks to make the deliveries. Their approach solves a deterministic model in each roll with the updated estimates of demand for the remaining periods. However, this deterministic model does not incorporate the evolving uncertain system states in terms of supply limits and demand. In this paper, we integrate the two-stage stochastic programming models within the rolling horizon approach to solve the post-disaster relief logistics planning problem. Our computational results demonstrate that this integration handles evolving stochastic system states more effectively.  

\subsection{Contributions to the Literature}
\label{sec:2.5}
Our contribution to the literature can be summarized as follows: 
\begin{itemize}
    \item Developing network flow type models for the last mile delivery problem.
    \item Modifying the deprivation cost calculation proposed by \citet{holguin_2013} to handle fluctuating demand over time.  
    \item Quantifying the value of stochastic programming by comparing it to the deterministic version of the problem. 
    \item Quantifying the value of the rolling horizon approach by comparing it to the solutions obtained from the static version of the solution approach.
\end{itemize}
In the remainder of the paper we will use RH to refer to the rolling horizon approach, S for the static approach, 2SSP to refer to the two-stage stochastic programming model, and D for the deterministic model. Thus, RH\_2SSP, RH\_D, S\_2SSP, and S\_D will represent the four different ways the post-hurricane logistics problem can be solved.

\section{Problem Setting}
\label{sec:3}
In this section, we first describe the logistics network on which the humanitarian relief operations planning problem is defined. Then the assumptions made in formulating the problem are described. Next, a deterministic optimization model and the corresponding notation are provided. Finally, we describe how various sources of uncertainty are modeled in the humanitarian relief logistics operations problem. 

\subsection{Logistics Network}
\label{sec:3.1}
\begin{figure}[htbp]
\begin{center}
  \includegraphics[scale=0.45]{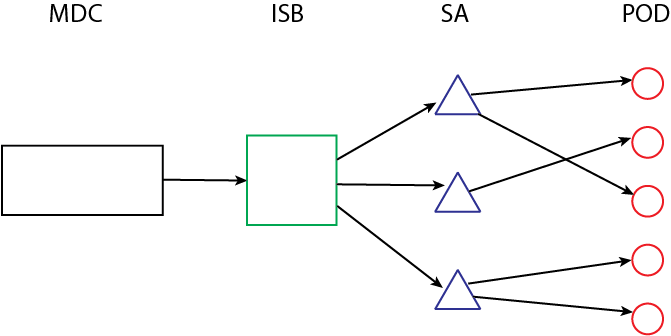}
\vspace{-0cm}
\caption{An illustration of the hierarchical structure of humanitarian relief logistics facilities.}
\label{fig:distribution_fa}       
\end{center}
\end{figure}

In the United States, large-scale disaster relief efforts are coordinated by FEMA. FEMA categorizes the disaster management lifecycle into four phases: (i) mitigation \--- taking actions to prevent or reduce the cause, impact, and consequence of future disasters, (ii) preparedness \--- making plans before an upcoming disaster, (iii) response \--- conducting operations to save lives and property immediately after a disaster, and (iv) recovery \--- rebuilding the communities back to normal  months or even years after a disaster. FEMA defines a hierarchical network structure of the distribution facilities involved in the preparedness and response phases of disaster logistics management. From the highest to the lowest level, this network consists of these components: Major Distribution Centers (MDCs), which are also referred to as logistics centers, followed by Incident Support Bases (ISBs), also referred to as pre-staging areas, then federal/state Staging Areas (SAs), and finally Points of Distribution (PODs) \citep{afshar_2012, vanajakumari_2016}. See Figure \ref{fig:distribution_fa} for an illustration of this network. 
\begin{itemize}
    \item The MDCs are permanent distribution facilities strategically located throughout a large geographical area used to store disaster relief commodities and recovery equipment. FEMA has a total of eight MDCs throughout the United States: Atlanta, GA; Ft. Worth, TX; Frederick, MD; Winchester, VA; Moffett Field, CA; Guam; Hawaii; and Caguas, Puerto Rico \citep{Ransom_2015}.  These distribution centers stockpile, among other items, millions of liters of water, millions of canned meals, and thousands of cots, blankets and tarps \citep{FEMA_2011}.
    \item When FEMA receives information about an advance-notice disaster, such as a flood or a hurricane, they will send some of the relief commodities forward from the MDCs to potential impact areas and set up the ISBs \citep{Smith_2011}. ISBs can enhance the ability of emergency management agencies to mobilize necessary supplies more quickly. 
    \item SAs are temporary facilities in the vicinity of areas projected to be affected by a disaster where commodities, equipment, and personnel are prepositioned for deployment \citep{afshar_2012}. Federal, state or county officials identify a set of potential SAs before the disaster. When an advance-notice disaster such as a hurricane is approaching, based on the predicted path, commodities are first forwarded from the MDCs to ISBs and then to SAs \citep{vanajakumari_2016}. 
    \item PODs are temporary local facilities within the affected  area and are operated by state officials. Relief commodities are distributed to the affected population through the PODs \citep{afshar_2012, vanajakumari_2016}.
\end{itemize}

\subsection{Assumptions}
\label{sec:3.2}
Before we formally present the proposed optimization models and their mathematical programming formulations, it is important to understand the key assumptions on which our work is based.
\begin{enumerate}
    \item The stochastic and deterministic models we propose address the last mile delivery problem in the humanitarian relief logistics operations for advance-notice disasters. We use hurricanes as an example through the rest of the paper. Network flow models are developed to solve the problem with uncertainty in demand and supply limits. Given our focus on the disaster relief network flow operations under demand and supply limits uncertainty, travel time and speed of distributing relief items are ignored for simplicity of modeling and analysis. 
    \item Of the four disaster management phases, we focus on the preparedness and response phases. More specifically, period $0$ is defined as the time period right before the hurricane strikes. The response phase is divided into $T$ time periods, which constitute the planning horizon of the proposed models. Transportation of relief items from the ISB to the SAs begins in period $0$, while transportation of relief items between SAs and PODs (and among themselves) begins in period $1$. They both continue throughout the planning horizon until period $T$. Each transshipment takes exactly one time period.
    \item Since our focus is on the preparedness and response phases, we consider a logistics network with one ISB, a set of PODs, and a candidate set of SAs. Determining when and where SAs should be opened is part of the decision variables. On the other hand, we assume that the PODs are ready for use immediately after a hurricane's landfall. Once an SA is opened, it remains open until the end of the planning horizon, and demand occurs at the PODs. In other words, we do not consider the transportation/delivery cost from PODs to individual beneficiaries.
    \item A delivery to a POD takes place only if the amount delivered plus the remaining inventory from the previous period can be used to fulfill the demand for at least one period. If no delivery is made to a POD then inventory is used to fully or partially fulfill the demand, depending on the amount of available inventory. We assume that partial fulfilment is recorded as demand not being satisfied. To handle such situations in our model, we assume that all remaining inventory at the PODs will head to a dummy node at no cost, if the inventory amount is insufficient to meet the demand for at least one period.
    \item Ground transportation is utilized for shipments from the ISB to the SAs and between SAs and PODs. Due to capacity constraints on ground transportation, emergency shipments utilizing air transportation may be needed during high demand periods from the ISB to the SAs or directly to the PODs. We have no capacity constraints for air transportation but it is a much more expensive mode of transportation.  
\end{enumerate}

\subsection{A Mixed Integer Programming Formulation}
\label{sec:3.3}
As mentioned in Section \ref{sec:2.5}, we develop two approaches (RH, S) and two models (2SSP, D) to solve the relief logistics problem. The 2SSP and D models are formulated as mixed integer programs (MIP). The deterministic (D) model is essentially a special case of the 2SSP model with only one scenario (\emph{i.e.}, the uncertain parameters are replaced with their corresponding mean values).

\subsubsection{Parameters and Decision Variables}
\label{sec:3.3.1}
The parameters used in the problem formulation are given in Table \ref{table:Indices}, and the decision variables are listed in Table \ref{table:two-stage SP_dv}. 

\begin{table}[htbp]\caption{Problem Parameters}\label{table:Indices}
\setlength{\tabcolsep}{1pt}
\footnotesize
\begin{tabular}{>{\raggedleft}p{0.05\textwidth}>{\centering}p{0.01\textwidth}>{\arraybackslash}p{0.94\textwidth}}
\hline\noalign{\smallskip}
$\mathcal{L}$&: & Set of potential SAs ($i,j \in\mathcal{L}=\{1,\ldots,L\}$) \\
$\mathcal{S}$&: & Set of available PODs ($i,j \in\mathcal{S}=\{L+1,\ldots,L+S\}$) \\
$\mathcal{K}$&: & Set of relief commodities ($k \in\mathcal{K}=\{1,\ldots,K\}$) \\
$\mathcal{T}$&: & Set of decision epochs ($t \in\mathcal{T}=\{0,\ldots,T\}$) \\
\noalign{\smallskip}\hline
$B_{ij}$&:  & Unit transportation cost of shipping relief items from SA or POD $i$ to SA or POD $j$ \\
$B^{g}_{i}$&:  & Unit transportation cost of shipping relief items to SA or POD $i$ from the ISB by ground transportation \\
$B^{h}_{i}$&:  & Unit transportation cost of shipping relief items to SA or POD $i$ from the ISB by air transportation \\
$\eta_{i}$&:  & Fixed cost of opening SA $i$ \\
$\zeta_{k}$&:  & Unit handling cost for relief item type $k$ \\
$\phi_{ik}$&:  & Capacity of SA or POD $i$ for relief item type $k$\\
$D_{ikt}$&: & Demand for relief item type $k$ in period $t$ at POD $i$\\
$R_{kt}$&: & Supply limit of relief item type $k$ in period $t$ from the ISB using ground transportation \\
\noalign{\smallskip}\hline
\end{tabular}
\end{table}

We separate our decision variables into state variables and local control variables. The state variables, treated as the first-stage decision variables in the 2SSP model, define the system status. The state variables ensure that we have an implementable decision for the actual post-hurricane under any realization of the stochastic process. We define $x$ variables as state variables to represent the status of the SAs (\emph{i.e.}, whether or not they are open), and $y$ variables as the actions of opening the SAs. In addition, state variables $\alpha$ and $U$ are defined to determine if demand is not met and for how long. Finally, state variable $V$ is used to describe the amount of the inventory for the SAs and PODs. The local control variables are the second-stage decision variables for the specific transportation plans among the ISB, the SAs, and the PODs. These decision variables are summarized in Table \ref{table:two-stage SP_dv}. Figure \ref{fig:Decision_Vars} illustrates these decision variables on a sample logistics network.
\begin{table}[htbp]
\caption{Decision Variables}\label{table:two-stage SP_dv}
\setlength{\tabcolsep}{1pt}
\footnotesize
\begin{tabular}{>{\raggedleft}p{0.05\textwidth}>{\centering}p{0.01\textwidth}>{\arraybackslash}p{0.94\textwidth}}
\hline\noalign{\smallskip}
\multicolumn{3}{p{\linewidth}}{First-stage (state) variables:} \\
$x_{it}$&:  & 1 if the status of SA $i$ is open at the beginning of period $t$; 0 otherwise \\
$y_{it}$&:  & 1 if we decide to open SA $i$ at the beginning of period $t$; 0 otherwise \\
$\alpha_{ikt}$&: & 1 if period $t$ demand of item $k$ at POD $i$ is satisfied; 0 otherwise \\
$V_{ikt}$&:  & amount of inventory of item $k$ at SA or POD $i$ at the end of period $t$ \\
$U_{ikt}$&:  & the number of periods since the last time demand was fully satisfied for item $k$ at POD $i$ in period $t$ \\
\noalign{\smallskip}\hline
\multicolumn{3}{p{\linewidth}}{Second-stage (local control) variables:} \\
$f_{ijkt}$&:  & amount of item $k$ shipped from SA or POD $i$ to SA or POD $j$ in period $t$ \\
$g_{ikt}$&:  & amount of item $k$ shipped to SA $i$ from the ISB by ground transportation in period $t$ \\
$h_{ikt}$&:  & amount of  item $k$ shipped to SA or POD $i$ from the ISB by air transportation in period $t$ \\
\noalign{\smallskip}\hline
\end{tabular}
\end{table}


\tikzstyle{ISB} = [draw=green, line width=0.5mm, regular polygon, regular polygon sides=4, fill=white!20,text width=2em, text centered, minimum height=2em]
\tikzstyle{SA} = [draw=blue, line width=0.5mm, regular polygon,regular polygon sides=3, text width=0.7cm, fill=white!20, inner sep=0pt, minimum height=1em]
\tikzstyle{POD} = [draw=red, line width=0.5mm, circle, fill=white!20,text width=4em, text badly centered, text width=1.2cm, inner sep=0pt, minimum height=2em]
\tikzstyle{Dummy} = [draw=red, line width=0.5mm, style=dashed, circle, fill=white!20,text width=4em, text badly centered, inner sep=0pt, minimum height=2em]
\tikzstyle{line} = [draw, -latex']
\begin{figure*}[ht]
\centering
    \begin{tikzpicture}[node distance = 2cm,auto]
        \node [SA, node distance = 4cm] (SA1) {\footnotesize{SA 1\\[1.5mm]$V_{1kt}$ }};
        \node [SA, below of=SA1, node distance = 2cm] (SA2) {\footnotesize{SA 2\\[1.5mm]$V_{2kt}$}};
        \node [SA, below of=SA2, node distance = 2cm] (SA3) {\footnotesize{SA 3\\[1.5mm]$V_{3kt}$}};
        \node [ISB, left of=SA2, node distance = 4.5cm] (ISB) {\footnotesize{ISB}};
        \node [POD, above right = 1cm and 2.5cm of SA2] (POD4) {\footnotesize{POD 4 \\[1.5mm]$V_{4kt}$}};
        \node [POD, above right = -2cm and 2.5cm of SA2] (POD5) {\footnotesize{POD 5 \\[1.5mm]$V_{5kt}$}};
        \node [Dummy, right = 4.5cm of SA2] (dummy) {\footnotesize{dummy}};
        \coordinate[above right of=POD4] (POD4_);
        \coordinate[below right of=POD5] (POD5_);
        \draw [->, black, line width=0.3mm, out=50,in=130] (ISB.north) to node[]{\textbf{$h_{4kt}$}} (POD4.north);
        \draw [->, black, line width=0.3mm, out=-50,in=-130] (ISB.south) to node[below]{\textbf{$h_{3kt}$}} (SA3.south);
        \draw [->, black, line width=0.3mm] (ISB.east) to node[]{\textbf{$g_{3kt}$}}(SA3.west);
        \draw [->, black, line width=0.3mm] (ISB.east) to node[]{\textbf{$g_{1kt}$}}(SA1.west);
        \draw [->, black, line width=0.3mm, out=-30, in=30] (SA1.east) to node[]{\textbf{$f_{13kt}$}}(SA3.north east);
        \draw [->, black, line width=0.3mm] (SA1.east) to node[]{\textbf{$f_{14kt}$}}(POD4.west);
        \draw [->, black, line width=0.3mm] (SA3.east) to node[below]{\textbf{$f_{35kt}$}}(POD5.west);
        \draw [->, black, line width=0.3mm] (POD4.south) to node[]{\textbf{$f_{45kt}$}}(POD5.north);
        \draw [->, black, line width=0.3mm, style=dashed] (POD4.east) to node[]{\textbf{$f_{4nkt}$}}(dummy.north);
        \draw [->, black, line width=0.3mm, style=dashed] (POD5.east) to node[]{\textbf{$f_{5nkt}$}}(dummy.south);
        \draw [->, black, line width=0.3mm] (POD4.north east) to node[]{\textbf{$D_{5kt}$}} (POD4_);
        \draw [->, black, line width=0.3mm] (POD5.south east) to node[]{\textbf{$D_{5kt}$}} (POD5_);
    \end{tikzpicture}
    \vspace{-4mm}
    \caption{An Illustration of the Decisions Variables on a Sample Logistics Network.}
    \label{fig:Decision_Vars}
\end{figure*}
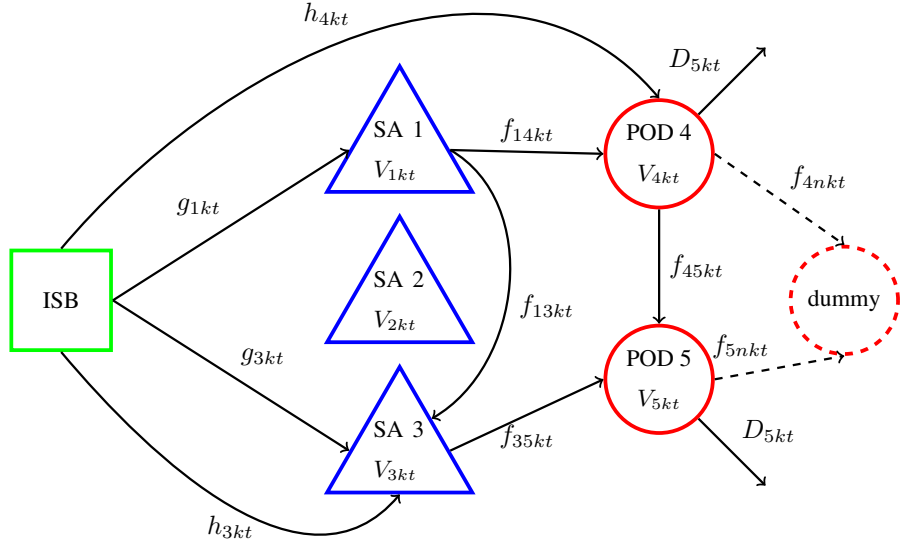

\subsubsection{Calculation of Deprivation Costs under Fluctuating Demand}
\label{sec:3.3.2}
As mentioned in Section~\ref{sec:2.2}, the deprivation cost is used to capture how human suffering increases as a function of the deprivation time. In their paper, \citet{perez_2016} simply computed the deprivation cost based on a function that only depends on the number of periods since the last period when the demand is fully satisfied, scaled by the demand amount. This function works well when demand is constant over time. However, as discussed in Section~\ref{sec:2.5}, the post-hurricane demand is likely to fluctuate over time due to the uncertainty in the network caused by damage to the infrastructure and inventories. Thus, we propose a new way to compute the deprivation cost that accommodates the fluctuating demand. 

Let $dep[i, k, U_{ikt}, t]$ denote the deprivation cost at POD $i$ when the demand for item $k$ goes unsatisfied for $U_{ikt}$ consecutive periods by time $t$. Each combination of $i, k, U_{ikt}$, and $t$ corresponds to a particular demand list, which can be used to calculate $dep[i, k, U_{ikt}, t]$. For example, if $T=5$, and the demand for relief item $k$ at POD $i$ is as given in Table \ref{table:exam_demand}, then the demand list for $[i, k, 3, 4]$ is $[120, 100, 120]$, which indicates that by period 4 the demand for relief item $k$ at POD $i$ was not satisfied in period 2, 3, and 4 but it was satisfied in period 1. 
\begin{table}[htbp]\caption{Example Demand }\label{table:exam_demand}
\centering
\setlength{\tabcolsep}{3pt}
\footnotesize
\begin{tabular}{ccccccc}
\noalign{\smallskip}\hline
&\multicolumn{6}{c}{Period \#} \\ \cline{2-7}
& period 0 & period 1 &period 2 &period 3 &period 4 &period 5   \\
\noalign{\smallskip}\hline\noalign{\smallskip}
Demand &0&100&120&100&120&130 \\
\noalign{\smallskip}\hline
\end{tabular}
\end{table}

Given a demand list for a relief item, Algorithm~\ref{alg:1} below can be used to calculate the deprivation cost under fluctuating demand. The algorithm implements a recursive function where $Dep$ is used to store the deprivation cost in the current stage, which is initially zero. In each recursive call of the algorithm, the minimum value in list $DL$ is recorded as $d_{min}$. In our example, $100$ is the minimum value in the demand list $[120, 100, 120]$, which indicates that $100$ units for relief item $k$ remain unsatisfied through periods 2\---4. Then, we update $Dep$ by adding the additional deprivation cost, which is equal to the product of $d_{min}$ and the deprivation cost function for relief item $k$, $\lambda_{k}(\tau)$, where $\tau$ is the length of the input demand list. Next, we update $DL$ by subtracting $d_{min}$ from each element in the list and save the new list into $DL^{'}$. Meanwhile, we capture the number of zeros among the elements in $DL^{'}$ and save it into $num\_{0}$. If all elements in $DL^{'}$ are zeros, we will return the current $Dep$ value and end the current recursive procedure; otherwise, we will separate $DL^{'}$ into sub-lists using the zeros as the break points and call function $compute\_dep$ recursively on these sub-lists. The recursive function will end after all unsatisfied demand has been fully examined. For the demand list in our example, which corresponds to $dep[i, k, 3, 4]$, the deprivation cost should be $100\times \lambda_{k}(3) + 2\times 20\times \lambda_{k}(1)$. 
\vspace{-0.4cm}
\begin{algorithm}[htbp]
	\caption{compute$\_$dep(demand list, $k$)}\label{alg:1}
	\begin{algorithmic}[gdg]
	    \State \textbf{Initialization:} Let $Dep \leftarrow 0$, $DL \leftarrow \text{demand list}$
	    \State $d_{min}\leftarrow \min(DL)$
	    \State $Dep \leftarrow Dep + d_{min}\times \lambda_{k}(len(DL))$
	    \State $DL^{'} \leftarrow DL - d_{min}$
	    \State $num\_{0} \leftarrow$ number of zeros in $DL^{'}$
	    \State $new\_list \leftarrow \emptyset$
	    \If {$num\_{0} = len(DL)$}
	    \State return $Dep$
	    \Else
	    \State Separate $DL^{'}$ into sub-lists using the zeros as the break points, save to $new\_list$
	    \For{each sub-list  $sl \in new\_list $}
		\State $Dep \leftarrow Dep + compute\_dep(sl, k)$
		\EndFor
		\State return $Dep$
	    \EndIf
	\end{algorithmic} 
\end{algorithm}
\vspace{-0.5cm}

\subsubsection{The First-stage Problem}
\label{sec:3.3.3}
We are now ready to present the first-stage problem. In the first-stage problem, we optimize over the state variables to minimize the fixed cost of opening SAs, handling cost of relief items, and the deprivation cost at PODs. The formulation for the first-stage problem (with first-stage cost only in the objective) is given as follows: 

\begin{subequations}\label{1st}
\begin{align}
    Min\quad z_{1}=&\sum\limits_{t\in \mathcal{T}} \Big( \sum\limits_{i\in \mathcal{L}}\eta_{i}y_{it} + \sum\limits_{i\in \mathcal{L}}\sum\limits_{k\in \mathcal{K}} \zeta_{k}V_{ikt} \Big) + \sum\limits_{i\in \mathcal{S}}\sum\limits_{k\in \mathcal{K}}dep[i, k, U_{ikT}, T] \nonumber\\
    &+ \sum\limits_{i\in \mathcal{S}}\sum\limits_{k\in \mathcal{K}}\sum\limits_{t\in \mathcal{T}\setminus\{0\}}dep[i, k, U_{ik(t-1)}, t-1]\alpha_{ikt}\label{1st-obj}\\
    V_{ik0}, U_{ik0}=& 0, \label{1st-ini} \qquad\qquad \forall i\in \mathcal{S},k \in \mathcal{K} \\
    x_{it}= &\sum^{t}_{t^{'} = 0}y_{it^{'}}, \label{1st-x-def} \qquad\qquad \forall i\in \mathcal{L}, t\in \mathcal{T}\\
    \sum_{t\in \mathcal{T}}y_{it} \leq &1, \label{1st-y-def} \qquad\qquad \forall i\in \mathcal{L}\\
    V_{ikt}\leq &\phi_{ik} x_{it}, \label{1st-SA-cap} \qquad\qquad \forall i\in \mathcal{L}, k \in \mathcal{K}, t\in \mathcal{T}\\
    U_{ikt}= &(1 - \alpha_{ikt})(U_{ik(t-1)}+1), \label{1st-U-def} \qquad\qquad \forall i \in \mathcal{S}, k \in \mathcal{K}, t\in \mathcal{T}\setminus\{0\}\\
    V_{ikt}\leq &\alpha_{ikt} \phi_{ik}, \label{1st-POD-cap-satis} \qquad\qquad \forall i \in \mathcal{S}, k \in \mathcal{K}, t\in \mathcal{T}\\
    x_{it},y_{it}\in& \{0,1\}, \label{1st-xy-domain} \qquad\qquad \forall i\in \mathcal{L}, t\in \mathcal{T}\\
    \alpha_{ikt}\in& \{0,1\}, \label{1st-alpha-domain} \qquad\qquad \forall i \in \mathcal{S}, k \in \mathcal{K}, t\in \mathcal{T}\\
    U_{ikt}\in& \mathbb{Z^{+}}, \label{1st-U-domain} \qquad\qquad \forall i \in \mathcal{S}, k \in \mathcal{K}, t\in \mathcal{T}\\
    V_{ikt}\geq& 0, \label{1st-V-domain} \qquad\qquad \forall i \in \mathcal{L} \cup \mathcal{S}, k \in \mathcal{K}, t\in \mathcal{T}
\end{align}
\end{subequations}

The objective function~\eqref{1st-obj} captures the total cost of the first-stage problem. The first part of the objective function corresponds to the facility cost, which includes the opening cost for the SAs and the inventory handling cost for the entire planning horizon. The second and third term compute the deprivation costs at the end of the planning horizon ($t=T$) and during the planning horizon, respectively. Next, we discuss the constraints:
\begin{itemize}
\item Constraints~\eqref{1st-ini} initialize the inventory levels at the PODs to zero and ensure that the demand for period $0$ is automatically satisfied. 
\item Constraints~\eqref{1st-x-def} ensure that the status of an SA is open if it was indeed opened at any point in time prior to or at the current period.
\item Constraints~\eqref{1st-y-def} make sure that an SA is opened at most once throughout the whole planning horizon so that the fixed cost is not incurred more than once.
\item Constraints~\eqref{1st-SA-cap} restrict the inventory of each relief item at an open SA to its capacity. 

\item Constraints~\eqref{1st-U-def} keep track of the number of periods that the demand goes unsatisfied for a relief item at a POD based on the following logic:

\begin{center}
         if $\alpha_{ikt} = 1$ then  $U_{ikt} = 0;$ otherwise  $U_{ikt} = U_{ik(t-1)} + 1$
\end{center}

\item Constraints~\eqref{1st-POD-cap-satis} ensures that $V_{ikt}  = 0$ if $\alpha_{ikt} = 0$. Recall that when the demand at a POD is not satisfied then any remaining inventory is ``artificially'' moved to a dummy node according to our assumption.
\end{itemize}

\subsubsection{The Second-stage Problem}
\label{sec:3.3.4}
The second-stage problem aims to minimize the total transportation cost of the relief effort given the first-stage decisions, demand, and the amount of supply limits. The second-stage problem is formulated as follows:
\begin{subequations}\label{2nd}
\begin{align}
    Min\quad z_{2} =& 
    \sum\limits_{k\in \mathcal{K}}\sum\limits_{t\in\mathcal{T}} \Big(\sum\limits_{i\in \mathcal{L}} B^{g}_{i} g_{ikt} + \sum\limits_{i\in \mathcal{L}\cup \mathcal{S}} B^{h}_{i} h_{ikt} + \sum\limits_{i\in \mathcal{L}\cup \mathcal{S}}\sum\limits_{j\in \mathcal{L}\cup \mathcal{S}, j\neq i} B_{ij} f_{ijkt}\Big) \label{2nd-obj} \\
    V_{ik0}= &g_{ik0}, \label{2nd-V-init} \qquad\qquad \forall\: i\in \mathcal{L}, k \in \mathcal{K}\\
    h_{ik0}= &0, \label{2nd-h-init} \qquad\qquad \forall\: i\in \mathcal{L}\cup \mathcal{S},k \in \mathcal{K}\\
    \sum\limits_{{\scriptscriptstyle j\in \mathcal{L} \cup \mathcal{S}, j \neq i }}f_{ijk0}= & 0, \label{2nd-f-init} \qquad\qquad \forall\: i \in \mathcal{L} \cup \mathcal{S}, k \in \mathcal{K}\\
    g_{ikt} + h_{ikt}\leq & V_{ikt}, \label{2nd-supply-cap} \qquad\qquad \forall\: i\in \mathcal{L}, k \in \mathcal{K}, t\in \mathcal{T}\\
    \sum\limits_{i\in \mathcal{L}} g_{ikt}\leq & R_{kt}, \label{2nd-gsupply-cap} \qquad\qquad \forall\: k \in \mathcal{K}, t\in \mathcal{T}\\
    V_{ikt}+ \sum\limits_{{\scriptscriptstyle j\in \mathcal{L} \cup \mathcal{S}, j \neq i}}f_{ijkt}= &V_{ik(t-1)}+ g_{ikt} + h_{ikt} \nonumber\\
    &+ \sum\limits_{{\scriptscriptstyle j\in \mathcal{L} \cup \mathcal{S}, j \neq i}}f_{jikt}, \label{2nd-SA-flow} \qquad\qquad \forall\: i \in \mathcal{L}, k \in \mathcal{K}, t\in \mathcal{T}\backslash\{0\}\\
    \sum\limits_{{\scriptscriptstyle j\in \mathcal{L} \cup \mathcal{S}, j \neq i }}f_{ijkt} \leq & V_{ikt}, \label{2nd-f-satis1} \qquad\qquad \forall\: i \in \mathcal{S}, k \in \mathcal{K}, t\in \mathcal{T}\backslash\{0\}\\
    \sum\limits_{{\scriptscriptstyle j\in \mathcal{L} \cup \mathcal{S}, j \neq i }}f_{jikt}\leq & V_{ikt}, \label{2nd-f-satis2} \qquad\qquad \forall\: i \in \mathcal{S}, k \in \mathcal{K}, t\in \mathcal{T}\backslash\{0\}\\
    h_{ikt}\leq & V_{ikt}, \label{2nd-h-satis} \qquad\qquad \forall\: i \in \mathcal{S}, k \in \mathcal{K}, t\in \mathcal{T}\backslash\{0\}\\
    V_{ikt} + \sum\limits_{{\scriptscriptstyle j\in \mathcal{L} \cup \mathcal{S} \cup \{n\}, j \neq i}}f_{ijkt}= &V_{ik(t-1)}- \alpha_{ikt}  D_{ikt} + h_{ikt} \nonumber\\
    &+ \sum\limits_{{\scriptscriptstyle j\in \mathcal{L} \cup \mathcal{S}, j \neq i}}f_{jikt}, \label{2nd-POD-flow} \qquad\qquad \forall\: i \in \mathcal{S}, k \in \mathcal{K}, t\in \mathcal{T}\backslash\{0\} \\
    g_{ikt}\geq& 0, \label{2nd-g-domain}  \qquad\qquad \forall\: i \in \mathcal{L}, k \in \mathcal{K}, t\in \mathcal{T} \\
    h_{ikt}\geq& 0, \label{2nd-h-domain} \qquad\qquad \forall\: i\in \mathcal{L}\cup \mathcal{S}, k \in \mathcal{K}, t\in \mathcal{T}\\
    f_{ijkt}\geq&0, \label{2nd-f-domain} \qquad\qquad \forall\: i, j\in \mathcal{L}\cup \mathcal{S} \cup \{n\}, k \in \mathcal{K}, t\in \mathcal{T}
\end{align}
\end{subequations}

The objective function~\eqref{2nd-obj} of the second-stage problem models the total transportation costs for humanitarian relief logistics operations, including ground and air transportation cost from the ISB to the SAs and between the SAs and the PODs. Next we discuss the constraints for the second-stage problem:
\begin{itemize}
\item Constraint sets~\eqref{2nd-V-init} and~\eqref{2nd-h-init} ensure that all inventory at the SAs in period $0$ is delivered by ground transportation and no air transportation is used for period $0$. Constraint set~\eqref{2nd-f-init} ensures that no delivery occurs between the SAs and the PODs in period 0.
\item Constraint set~\eqref{2nd-supply-cap} ensures that the amount of relief items delivered by ground and air transportation does not exceed the available inventory level. As mentioned in Section~\ref{sec:3.2}, the amount shipped via ground transportation should not exceed the supply limits as defined in the constraint set~\eqref{2nd-gsupply-cap}. The flow balance constraints for the SAs are given by \eqref{2nd-SA-flow}. 
\item Constraint sets~\eqref{2nd-f-satis1}-\eqref{2nd-h-satis} ensure that shipments to and from PODs are bounded by the available inventory level, and the flow balance constraints for the PODs are given by \eqref{2nd-POD-flow}. 
\end{itemize}

\subsection{Two-stage Stochastic Programming Models}
\label{sec:3.4}

\subsubsection{Stochastic Model \iffalse(Markov chain, and stagewise independence given the Markovian states, etc.)\fi}
\label{sec:3.4.1}
As mentioned in Section~\ref{sec:2}, demand from affected populations can fluctuate due to a variety of uncertain factors caused by damage from the hurricane. The uncertainty in supply availability can stem from delays or losses of relief goods in the supply chain. To generate more realistic demand and supply limits, we use a Markov chain to generate input data for our stochastic model. The idea of applying the Markov chain model is inspired by~\citet{regnier_2006}, who studied public evacuation problems using information given by a Markov chain that forecasts the evolution of a hurricane. Also, \citet{taskin_2010} used a Markov chain to predict hurricane landfall counts for an upcoming season and developed a multi-period stochastic programming model to guide production planning for a manufacturer. Our process of generating input data contains two layers. First, we generate a sequence of Markovian states, each of which represents the damage level for future periods based on the underlying Markov chain model. Then, for each observed state on the sequence, random demand and supply limits are generated by two independent probability distributions conditioned on the observed state. The details are presented below.
\begin{itemize}
\item Layer 1: We use a sequence of Markovian states to describe possible damage levels to the relief logistics network at each decision epoch during the post-hurricane response phase. The state space of the Markov Chain $\mathcal{M}$ includes all possible damage levels of the network. For example, we can define three damage levels: high(H), medium(M), and low(L) as the state space $\mathcal{M} =$ \{H, \:\:M, \:\:L\} of the Markov chain. Also, $\mathcal{P}$ is used to store the one-step transition probability matrix of the Markov chain. 
\item Layer 2: The second layer is the realization of supply limits and demand values based on a conditional probability distribution given the specific Markovian state. We use two independent probability distributions conditional on each Markovian state $m\in \mathcal{M}$, for the associated demand and supply limits. Also, the realization process has stagewise independence, which means the realization at the current point on the Markovian states sequence does not affect the realization for remaining points on the sequence.

\end{itemize}
Applying this two-layer stochastic process, we can generate a scenario tree for use in our proposed 2SSP models, which we present next. 

\subsubsection{Static 2SSP Formulation}
\label{sec:3.4.2}
In Section \ref{sec:3.3} we provided the formulation for the deterministic version of our problem which is a special case of the 2SSP model. Here, we describe how we generalize the deterministic formulation that leads to the 2SSP formulation. First, let $\Omega$ be the set of all possible scenarios and $P^{\omega}$ the probability of scenario $\omega \in \Omega$. Then $D^{\omega}_{ikt}$, $R^{\omega}_{kt}$, $f^{\omega}_{ijkt}$, $g^{\omega}_{ikt}$, and $h^{\omega}_{ikt}$ are the same variables as defined in Table \ref{table:two-stage SP_dv} but under scenario $\omega$.

The objective function~\eqref{static-obj} of the static 2SSP model aims to minimize the total expected cost, including the transportation, inventory, and deprivation costs as formulated below.

\begin{subequations}\label{static}
\begin{align}
    Min \quad z = &\sum\limits_{t\in \mathcal{T}} \Big(\sum\limits_{i\in \mathcal{L}}\eta_{i}y_{it} + \sum\limits_{i\in \mathcal{L}}\sum\limits_{k\in \mathcal{K}} \zeta_{k}V_{ikt}\Big) \nonumber\\
    &+ \sum\limits_{\omega \in \Omega} P^{\omega} \sum\limits_{k\in \mathcal{K}}\sum\limits_{t\in \mathcal{T}} \Big(\sum\limits_{i\in \mathcal{L}} B^{g}_{i} g^{\omega}_{ikt} + \sum\limits_{i\in \mathcal{L}\cup \mathcal{S}} B^{h}_{i} h^{\omega}_{ikt} + \sum\limits_{i\in \mathcal{L}\cup \mathcal{S}}\sum\limits_{j\in \mathcal{L}\cup \mathcal{S}, j\neq i} B_{ij} f^{\omega}_{ijkt}\Big)\nonumber\\
    &+ \sum\limits_{\omega \in \Omega} P^{\omega} \Big(\sum\limits_{i\in \mathcal{S}}\sum\limits_{k\in \mathcal{K}}\sum\limits_{t\in \mathcal{T}}dep[w, i, k, U_{ik(t-1)}, t-1]\alpha_{ikt}\nonumber\\
    &+ \sum\limits_{i\in S}\sum\limits_{k\in K}dep[w, i, k, U_{ikT}, T]\Big)\label{static-obj}
\end{align}
\end{subequations}
The first-stage problem for 2SSP is same as the deterministic first-stage problem provided in Section \ref{sec:3.3.3}. The only difference is that the objective function \eqref{1st-obj} is replaced with \eqref{static-obj}. However, the second-stage problem needs to be revised to be a scenario-based problem. Thus, every local control variable in the static second-stage problem is associated with the scenario index $\omega$. The full formulation for the static 2SSP formulation can be found in the Appendix. 

\subsubsection{Rolling Horizon 2SSP Approach}
\label{sec:3.4.3}
\begin{figure}[htbp]
\begin{center}
  \includegraphics[scale=0.45]{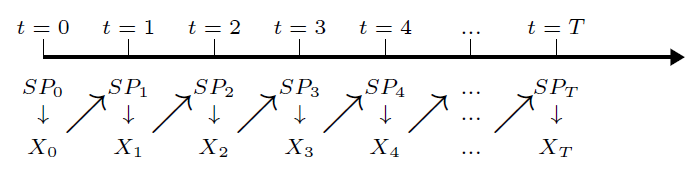}
\vspace{-0cm}
\caption{An Illustration of the Rolling Horizon 2SSP Approach.}
\label{fig:2}       
\end{center}
\end{figure}

One potential drawback of the static 2SSP formulation is that it is inadequate to address the dynamic evolution of the aforementioned uncertain factors in the problem. To this end, we propose a stochastic lookahead framework that implements the 2SSP model in a rolling horizon fashion. Specifically, as shown in Figure~\ref{fig:2}, the rolling horizon approach separates the timeline into sub-problems for each period ($SP_{t}$). For each sub-problem in period $t < T$, we run a static 2SSP model in which the state variables define the logistics system states from period $t$ to period $T$, and the second-stage decisions are the humanitarian relief logistics operation plans for the remaining periods in the planning horizon. After solving the model, only the logistics plans for the current period, $X_t$, are implemented, and subsequently the horizon is shifted. At the beginning of the next period, the state of the systems is observed and a new optimization problem is solved. The actual hurricane relief logistics operations for the entire time horizon will then be \{$X_{0}, X_{1}, X_{2},..., X_{T}$\}, which are used to compute the total cost of the rolling horizon 2SSP approach. Compared with the rolling horizon approach using deterministic optimization models~\citep{rivera_2016, vanajakumari_2016, hasani_2018} mentioned in Section~\ref{sec:2.4}, our proposed rolling horizon 2SSP approach better handles the evolving uncertainty during the post-hurricane response period because it considers a more comprehensive future trajectory. The details of the comparisons of RH\_2SSP, RH\_D, S\_2SSP, and S\_D are provided in our numerical experiment results shown in Section~\ref{sec:4.3}. 

\section{Numerical Experiments}
\label{sec:4}
We conducted extensive numerical experiments to evaluate the performances of the proposed static and rolling horizon approaches. We used a scenario tree to describe the stochastic process for the post-hurricane logistics network. The performance of each proposed approach is evaluated based on the respective total costs (including logistics and social costs) via an out-of-sample test. The details of the out-of-sample test are described in Section \ref{sec:4.3}.  We used Gurobi 7.0.2 to solve all the optimization models on the Clemson Palmetto Cluster\footnote{Palmetto Cluster is Clemson University’s primary high-performance computing (HPC) resource, which is utilized by researchers, students, faculty, and staff from a broad range of disciplines.}, where each compute node is equipped with Intel Xeon Processor E5410 8 cores @2.33 GHz and 16GB RAM.

\subsection{Problem Data}
\label{sec:4.1}
\begin{figure}[htbp]
\begin{center}
  \includegraphics[scale=0.42]{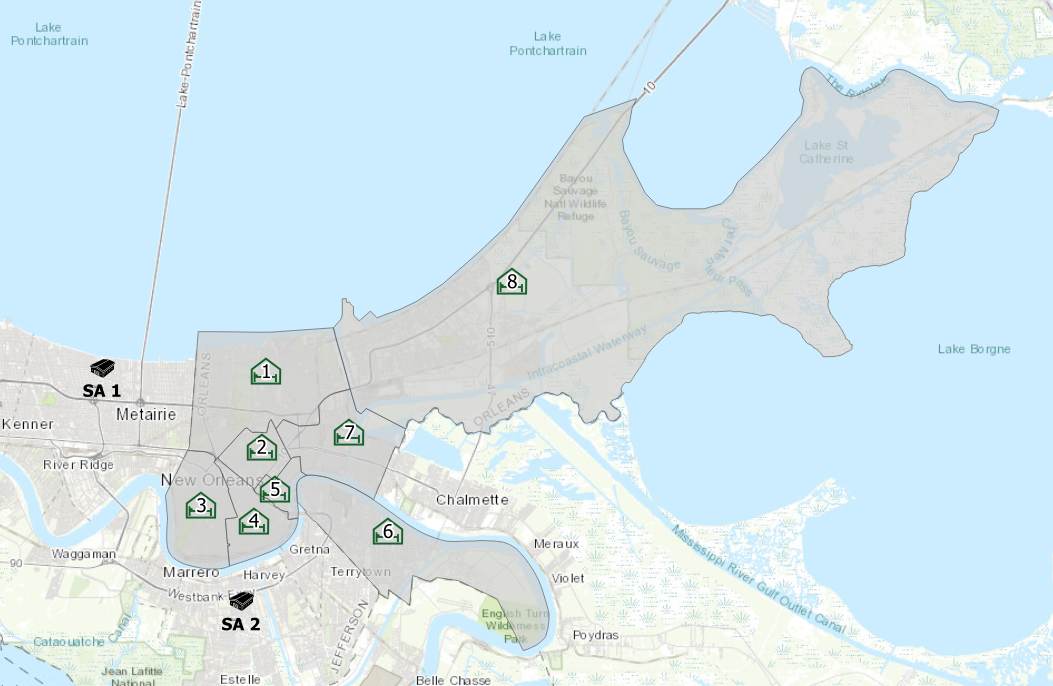}
\caption{Map of the Staging Areas and Points of Distributions.}
\label{fig:NO_map}       
\end{center}
\end{figure}

In our numerical experiments, we construct problem data based on the transportation network of City of New Orleans, which was used by \citet{perez_2016} in their study. As shown in Figure \ref{fig:NO_map}, the area is divided into eight districts by the New Orleans Police Department. We assume that each district has a POD located near its center which is ready to be used in the aftermath of a hurricane, so that we have eight PODs in total. Also, two potential SAs in Jefferson Parish are pointed out in the map: one located in Metairie and the other in Marrero. The length of the planning horizon for the response phase is set to be five periods ($T = 5$). The unit cost of shipping via ground transportation between SAs and PODs are proportional to the travel distances. Also, we consider two types of relief commodities in our experiments. These two commodities have the same deprivation cost function but different amounts of initial demand. 

Table~\ref{tab:para_val} summarizes the parameters used in our experiments. The deprivation cost function $\lambda_{k}(t)$ follows the structure given by \citet{perez_2016}. We introduce a weight parameter $\delta$ to represent the relative weight of the deprivation cost in the overall cost. As mentioned earlier, the ground transportation costs between different nodes are estimated based on their travel distances. Air transportation, which is significantly more expensive, is used as an alternative mode if needed. The unit air transportation cost is $50$ times that of ground transportation. This value is selected based on literature related to air transportation \citep{barbarosouglu_2002, rawls_2010, lin_2012, chowdhury_2017}. The capacity of each POD for each relief item is 350, which is slightly larger than 1.5 times the largest baseline demand among all PODs shown in Table \ref{tab:base_demand}. The capacities of the SAs are considered to be equal to the sum of capacities of all eight PODs. The fixed cost for opening an SA and the handling cost for relief items are based on data from \citet{vanajakumari_2016}.
\begin{table}[htbp]
\caption{Parameters Used in Numerical Experiments}
\label{tab:para_val}       
\centering
\begin{tabular}{lc}
\hline\noalign{\smallskip}
Deprivation Cost Function ($\lambda_{k}$) & $\lambda_{k}(t) = \delta \cdot (e^{1.5+0.12\cdot t}-e^{1.5})$\\
Deprivation Cost Parameter ($\delta$)  & 30\\
Facility capacity for each relief item ($\phi$)& \\
\quad SA & 1400\\
\quad POD & 350 \\
Unit delivery cost of relief items &\\
\quad  from ISB to SAs by ground ($B^{g}$) & 10 \\
\quad from ISB to SAs/PODs by air ($B^{h}$) & 500\\
Unit handling cost of relief items ($\zeta$)& 1\\
Fixed cost for opening an SA ($\eta$) & 10,000\\
\noalign{\smallskip}\hline
\end{tabular}\\
\end{table}\\

\subsection{Scenario Tree Generation}
\label{sec:4.2}
We use a scenario tree representation to capture the uncertainty in demand for and supply limits of relief items following a hurricane. We consider a two-layer stochastic model. As mentioned in Section~\ref{sec:3.4.1}, the first layer is the Markovian states, and the second layer is the realization of the random parameters based on a probability distribution conditional on each Markovian state. 

\begin{figure}[htbp]
\begin{center}
  \includegraphics[scale=0.35]{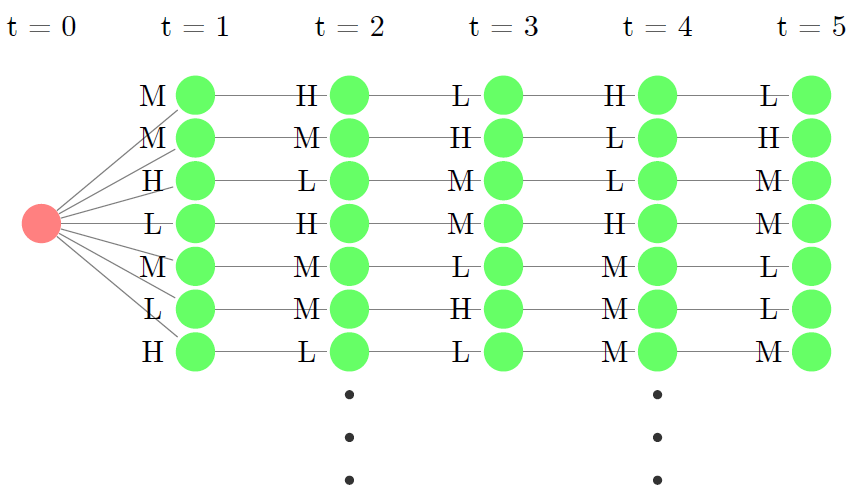}
\vspace{-0cm}
\caption{An Illustration of a Set of Sample Paths}
\label{fig:4}       
\end{center}
\end{figure}

In the first layer, we assume that during each decision epoch the logistics system is in one of three states: H, M, and L. In state H, the demand at each POD is more likely to be higher than the baseline value and supply limits more likely to be lower; in state L, the demand is likely to be lower and supply limits higher than the corresponding baseline values; in state M, both demand and supply limits are close to the baseline values. Table \ref{tab:base_demand} shows the baseline demand values. These values were determined based on the population of each POD district. The baseline values for the supply limits were set to be the total baseline demand, \textit{i.e.}, the baseline supply limit is 1395 for relief item 1 and 1224 for item 2.
\vspace{-0.5cm}
\begin{table}[htbp]
\caption{Baseline Demand}
\label{tab:base_demand}
\begin{center}
    \begin{tabular}{ccccccccc}
     \hline
     {}&POD 1& POD 2 &POD 3 & POD 4 & POD 5 & POD 6 & POD 7 & POD 8\\
     \hline
     Item 1 &219 & 214& 195& 116& 70& 162& 205& 214\\
     Item 2 &177 & 163& 171& 121& 69& 130& 190& 203\\
     \hline
    \end{tabular}\\
\end{center}
\end{table}
\vspace{-0.5cm}
As discussed in Section~\ref{sec:3.4.1}, we model the post-hurricane humanitarian relief logistics network status over time as a stochastic process described by a Markov chain. We use synthetic data in our study to demonstrate the effectiveness of the proposed models and solution approaches, and we perform an extensive sensitivity analysis on different model parameters. In practice, estimating these problem parameters can be very challenging, but it is possible by utilizing historical data according to existing works in the literature. For example, the transition probability matrix of the Markov chain can be estimated based on historical data on past hurricanes' trajectories and their impacts to the affected area's logistics network. \citet{taskin_2010} used a discrete-time Markov chain model using 50 years of historical hurricane records to predict the demand scenarios for hurricane season and addressed a stochastic inventory control problem for manufacturing and retail firms. In our experiments, we assume that the system is initially in state H with probability 0.3, in state M with probability 0.4, and in state L with probability 0.3. The one-step transition probability matrix is shown below. Based on this transition probability matrix, for example, if the system is in state H at some time $t$ then at time $t+1$ the system will be in state H, M, or L with probabilities 0.7, 0.2, and 0.1, respectively. Figure \ref{fig:4} shows a set of sample paths generated according to the Markov chain model described above, where sample path corresponds to a sequence of Markovian states that would occur in the planning horizon. 
\begin{center}
    \[
\def\stackalignment{r}
\textbf{{$\mathcal{P}$}} = 
       \Centerstack{
        H  \\
        M  \\
        L    
        }\!
        \stackon{
        \stretchleftright{\bigg[\:}{\tabbedCenterstack{
         0.7 & 0.2 & 0.1 \\
        0.3 & 0.5 & 0.2 \\
         0.1   & 0.2   & 0.7   
        }}{\bigg]}}{\tabbedCenterstack{
        \protect\phantom{0.}H & \protect\phantom{0.}M & \protect\phantom{0.}L}
        \kern1pt}
\]
\end{center}


The second layer models the actual realizations of the demand and supply limits given the corresponding Markovian state. The distributions of supply limits and demand conditional on the Markovian state of the system are given in Tables \ref{tab:pro_supply} and \ref{tab:pro_d}, respectively. For example, if the system is in state H then the realization of the supply limits is equal to the baseline level (100\%) with probability 0.25, $110\%$ with probability 0.1, $120\%$ with probability 0.05, $90\%$ with probability 0.35, and $80\%$ with probability 0.25. The realization of the demand level is slightly different, it is assumed to follow a uniform distribution between two parameters that are determined following a probability distribution conditional on the Markovian states given in Table~\ref{tab:pro_d}. For example, if we are in state H then with probability 0.15 the demand at each POD is given by the corresponding baseline demand multiplied by a factor that is uniformly distributed between 1.35 and 1.45 (see the number in the upper left corner of Table \ref{tab:pro_d}). We assume that the stochastic demand and supply limits are independent random variables given a Markovian state. 
\begin{table}[htbp]
\caption{Probability Distributions for Supply Limit Conditioned on Markovian States\iffalse, $\mathcal{F'}$ {\color{blue} XXX Are $\mathcal{F'}$ and $\mathcal{F}$ ever used later? Why do you introduce these notations?? XXX}\fi}
\label{tab:pro_supply}
\begin{center}
    \begin{tabular}{p{1cm}p{1cm}p{1cm}p{1cm}p{1cm}p{1cm}}
     \hline
     & 120\%&110\%& 100\% &90\%&80\%\\
     \hline
     H & 0.05 & 0.10 & 0.25&0.35 & 0.25 \\
     M & 0.20 & 0.15 & 0.35 &0.20 & 0.10 \\
     L & 0.30 & 0.30 &0.25 & 0.10 & 0.05 \\
     \hline
    \end{tabular}\\
\end{center}
\end{table}\begin{table}[htbp]
\caption{Probability Distributions for Demand Levels Conditioned on Markovian States \iffalse, $\mathcal{F}$\fi}
\label{tab:pro_d}
\begin{center}
    \begin{tabular}{p{0.4cm}p{0.8cm}p{0.8cm}p{0.8cm}p{0.8cm}p{0.8cm}p{0.8cm}p{0.8cm}p{0.8cm}p{0.8cm}}
     \hline
     & +45\%&+35\%&+25\%&+15\%&+5\% &-5\%&-15\% &-25\%&-35\%\\
     & \~+35\%&\~+25\%&\~+15\%&\~+5\%&\~ -5\% &\~-15\%&\~-25\% &\~-35\%&\~-45\%\\
     \hline
     H & 0.15 & 0.3 & 0.25 & 0.15 & 0.05&0.04 & 0.04 & 0.02 & 0\\
     M & 0.04 & 0.06 & 0.1 & 0.25 & 0.25&0.1 & 0.1 & 0.06 & 0.04\\
     L & 0 & 0.02 & 0.04 & 0.04&0.1 & 0.15 & 0.2 & 0.3& 0.15\\
     \hline
    \end{tabular}\\
\end{center}
\end{table}

\subsection{Experimental Setup}
\label{sec:4.3}
To test the out-of-sample performance of the proposed static and rolling horizon approaches, we generated 1000 sample paths following the two-layer process we described in Section~\ref{sec:4.2}. The sample size of 1000 in the out-of-sample experiment is selected to capture enough variability for the stochastic model under an acceptable computational budget. For S\_2SSP, we first randomly sample 100 scenarios and solve the S\_2SSP model. Then we evaluate the obtained solution on the 1000 sample paths and record the performance metrics of interest. As will be shown in Section~\ref{sec:4.4.1}, we perform a stability test to justify that the sample size of 100 is sufficient. In the out-of-sample test for the RH\_2SSP, we use the same 1000 sample paths used for testing S\_2SSP. For each roll of the RH\_2SSP, 10 scenarios are randomly sampled based on the network state for the corresponding period. We justify that the sample size 10 is appropriate for our experiment, again by performing a stability test, which is described in Section~\ref{sec:4.4.2}. Note that every time we finish a roll, we update the state variables for the current and subsequent periods. For the next roll, a new set of 10 scenarios is randomly generated for the remaining periods based on the state of the system at that point in time. The computational time of the S\_2SSP(S\_D) is recorded as the computational time of a single S\_2SSP model with 100(1) scenarios, while for the RH\_2SSP and the RH\_D the computational time is recorded as the sum of computational time spent in each roll for the entire sample path, and the average computational time over all sample paths considered in the out-of-sample test is reported. To give reliable statistics, for each experiment we perform ten replications, and the results provided in the following sections correspond to the average results over these ten replications. To ease our exposition, we will use the abbreviations listed in Table \ref{tab:results_term} for the remainder of this section.

\begin{table}[htbp]
\caption{Abbreviations for Result Illustrations}
\label{tab:results_term}       
\centering
\footnotesize
\begin{tabular}{p{0.2\textwidth}p{0.7\textwidth}}
\hline\noalign{\smallskip}
A & = \{RH\_2SSP, S\_2SSP, RH\_D, S\_D\}\\
a\_Total & Total cost for a $\in$ A\\
a\_Dep & Deprivation cost for a $\in$ A\\
a\_Log & Logistics cost for a $\in$ A\\
a\_Air & Air transportation cost for a $\in$ A\\
a\_compT & Computational time for a $\in$ A\\
a\_Dep/Total & Ratio of deprivation cost to the total cost for a $\in$ A \\
Improv\_P & Percentage of improvement for RH\_2SSP\_Total compared to S\_2SSP\_Total \\
\noalign{\smallskip}\hline
\end{tabular}
\vspace{-0.5cm}
\end{table}

\subsection{Stability Test}
\label{sec:4.4}
\subsubsection{Stability Test for S\_2SSP}
\label{sec:4.4.1}
The stability test intends to show that using 100 scenarios for solving 2SSP models in the static approach is a reasonable sample size. To perform this stability test, we consider a set of sample sizes: 20, 50, 100, 200, and 400, and then evaluate the corresponding solutions through an out-of-sample test with 1,000 scenarios. The results summarized in Table~\ref{tab:var_diffn_sta} indicate that different sample sizes considered in 2SSP models lead to different total costs. Clearly, using larger sample sizes will get us closer to the true cost value. In addition, the standard deviation of the cost reduces as the sample size increases. However, using large sample sizes increases the computational time. Table~\ref{tab:var_diffn_sta} shows the results over ten replications. As can be seen in the table, the average cost stabilizes at about 100 scenarios. Using 200 scenarios reduces the standard deviation by about 50\% but the average cost changes only marginally and the computational time more than triples. Thus, in the results that we will present in the following sections we use 100 scenarios when solving the 2SSP models. 
\begin{table}[htbp]
\caption{The Performance of S\_2SSP using Different Sample Sizes for Solving 2SSP}
\label{tab:var_diffn_sta}       
\centering
\begin{tabular}{lrrrrr}
\hline\noalign{\smallskip}
& \multicolumn{5}{c}{Number of Scenarios} \\ \cline{2-6}
& \multicolumn{1}{c}{20} & \multicolumn{1}{c}{50} & \multicolumn{1}{c}{100} & \multicolumn{1}{c}{200} & \multicolumn{1}{c}{400} \\
\hline
S\_2SSP\_Total: Average (\$) &247,226 &224,472 &221,135 &220,145 & 219,985  \\
S\_2SSP\_Total: Standard Deviation (\$) &30,978 &8,186 &3,793 &1,532 &1,287  \\
S\_2SSP\_compT (seconds) &123 &679 &2,716 &8,987 &  30,261\\
\noalign{\smallskip}\hline
\end{tabular}
\end{table}

\subsubsection{Stability Test for RH\_2SSP}
\label{sec:4.4.2}
This stability test intends to show that using 10 scenarios for solving 2SSP models in the rolling-horizon approach is a reasonable sample size. To perform this stability test, we consider a set of sample sizes: 10, 20, 50, and 100, and evaluate the performances of the respective solutions using out-of-sample tests. As can be seen from Table \ref{tab:var_diffn_roll}, the overall costs do not vary much when different sample sizes are used. On the other hand, the computational time of the RH\_2SSP significantly increases as the sample size increases from 10 to 100. In summary, the solution quality is not sensitive to the number of scenarios used in each roll, whereas the computational time heavily depends on this sample size. Thus, in our numerical experiments we use 10 sample paths in each roll of the RH\_2SSP. 
\begin{table}[htbp]
\caption{The Performance of RH\_2SSP Using Different Sample Sizes for Solving 2SSP in Each Roll}
\label{tab:var_diffn_roll}       
\centering
\begin{tabular}{lrrrr}
\hline\noalign{\smallskip}
& \multicolumn{4}{c}{Number of Sample Paths} \\ \cline{2-5}
& \multicolumn{1}{c}{10} & \multicolumn{1}{c}{20} & \multicolumn{1}{c}{50} & \multicolumn{1}{c}{100} \\
\hline
RH\_2SSP\_Total: Average (\$)&195,564  &195,258 &195,561 &195,729   \\
RH\_2SSP\_Total: Standard Deviation (\$) &5,706 &4,499 &3,523 &3,238   \\
RH\_2SSP\_compT (seconds) &85 &1,165 &4,171 &7,406   \\
\noalign{\smallskip}\hline
\end{tabular}
\end{table}

\subsection{Experiment Results}
\label{sec:4.5}

\subsubsection{Value of Stochastic Programming}
\label{sec:4.5.1}
To show the value of a two-stage stochastic programming model compared to a deterministic model, we solved the deterministic and the 2SSP versions for both the static and the rolling horizon approaches. The demand and supply limits that we consider in the deterministic models are set to be the average values of the respective realizations of random variables used in their 2SSP model counterpart, \textit{i.e.}, average over the 100(10) scenarios for the S\_2SSP(RH\_2SSP). We compared the performance of static approach under the 2SSP and deterministic methods with different weights for the deprivation cost, and collected the results in Table \ref{tab:DeterStochas_1}. From this table we clearly see that the deterministic model can be solved very quickly but provides poor results. Specifically, the total cost for the deterministic model is almost 4.5 times higher compared to the 2SSP model when weight was 1, although it is significantly faster to solve. Even when we set the weight to 5 or 10, a large gap remains between the overall costs for the 2SSP and deterministic models. This is expected because the deterministic model considers only a single scenario representing the average value of random realizations, which fails to provide a robust decision for the state variables to accommodate fluctuating demand and supply limits during the out-of-sample test. Instead, it resorts to air transportation mode to meet the fluctuating demand, which incurs a much higher cost. 

\begin{table}[htbp]
\caption{Comparison of Using Deterministic Optimization Models and Stochastic Optimization Models for the Static Approach}
\label{tab:DeterStochas_1}       
\centering
\begin{tabular}{lrrrrrr}
\hline
&\multicolumn{6}{c}{Weight/Model} \\ \cline{2-7}
{}&\multicolumn{2}{c}{1}&\multicolumn{2}{c}{5}&\multicolumn{2}{c}{10} \\
&\multicolumn{1}{c}{S\_2SSP} &\multicolumn{1}{c}{S\_D} &\multicolumn{1}{c}{S\_2SSP} &\multicolumn{1}{c}{S\_D}&\multicolumn{1}{c}{S\_2SSP} &\multicolumn{1}{c}{S\_D} \\
\hline\noalign{\smallskip}
\_Total (\$) &221,135 &982,132 &532,754 &1,154,508 &676,648 &1,136,203  \\
\_Log (\$) &127,493 &960,397 &325,170 &1,154,508  &498,192 &1,136,203   \\
\_Dep (\$) &93,642 &21,735 &207,584 &0 &178,456 &0 \\
\_Air (\$) &2,615 &804,907 &151,184 &982,255 &309,654 &963,323   \\
\_Dep/Total &42.35\% & 2.13\% & 38.96\% &  0\%& 26.37\% &  0\%\\
\_compT (seconds) &2,716 &2 &59  &1  &42  &1   \\
\noalign{\smallskip}\hline
\end{tabular}
\end{table}

\begin{figure}[htbp]
\begin{center}
  \includegraphics[scale=0.385]{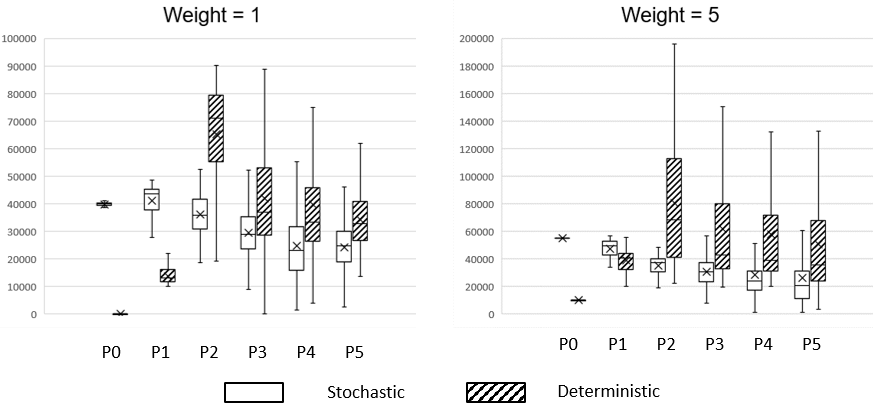}
\caption{Illustration of Costs per Period for the Rolling Horizon Stochastic and Deterministic Models}
\label{fig:boxplot}       
\end{center}
\end{figure}

We also compared the performance of rolling horizon approach under the 2SSP and deterministic methods under different weights for the deprivation cost. Table~\ref{tab:DeterStochas_2} summarizes the results where ``\_Diff" shows the cost reduction achieved by the RH\_2SSP compared to the RH\_D, which is calculated by $1-(\text{RH\_2SSP\_Total}/\text{RH\_D\_Total})$. We can see that the cost savings obtained by RH\_2SSP compared to RH\_D can be different under different weights. To provide more detailed information on what happens in each roll using the two approaches, we illustrate the cost distribution in each period using boxplots in Figure~\ref{fig:boxplot}. The boxplots on the left(right)-hand side indicate the cost for each roll under the 2SSP(D) model. For each boxplot, the cross mark represents the mean value, the central mark in the box is the median, and the upper and lower bound of each box correspond to the 25th and 75th percentiles, respectively. 

\begin{table}[htbp] 
\caption{Comparison of the Deterministic and Stochastic Models Using the Rolling Horizon Approach under Different Deprivation Cost Weights}
\label{tab:DeterStochas_2}       
\centering
\begin{tabular}{lrrrrrr}
\hline
&\multicolumn{6}{c}{Weight/Model} \\ \cline{2-7}
{}&\multicolumn{2}{c}{1}&\multicolumn{2}{c}{5}&\multicolumn{2}{c}{10} \\
&\multicolumn{1}{c}{RH\_2SSP} &\multicolumn{1}{c}{RH\_D}& \multicolumn{1}{c}{RH\_2SSP} &\multicolumn{1}{c}{RH\_D}& \multicolumn{1}{c}{RH\_2SSP} &\multicolumn{1}{c}{RH\_D}\\
\hline\noalign{\smallskip}
\_Total (\$) &195,564 &194,741 & 223,592& 300,111 &236,446 &419,702    \\
\_Log (\$) &174,299 &123,226 & 211,616& 180,807 &212,725 &184,289    \\
\_Dep (\$)&21,265 &75,515 &11,976 &119,304 &23,721 &235,413    \\
\_Air (\$) &1 &8 &196 &696 &955 &3,341    \\
\_Dep/Total &10.31\% & 38.78\%& 5.36\% &39.74\% & 10.03\%& 56.09\%\\
\_compT (seconds) &85 &1 &108 &4 &35 &1  \\
\_Diff  &\multicolumn{2}{c}{-0.42\%} &\multicolumn{2}{c}{25.49\%} &\multicolumn{2}{c}{43.66\%}   \\
\noalign{\smallskip}\hline
\end{tabular}
\end{table}

From the boxplot with weight equal to 1, we can observe that for the first two periods, \emph{i.e.}, periods 0 and 1, the mean costs of RH\_2SSP are much higher than that of RH\_D. Starting from period 2,  RH\_D results in a higher mean cost and the corresponding variance is also higher. Similar observations are seen when the weight parameter is set to 5. The boxplots of weight equal to 10 are omitted in Figure 6 because the output was very close to the one with weight equal to 5. When the weight is equal to 1, the total costs for RH\_D and RH\_2SSP are almost the same (RH\_D was 0.42\% lower), but RH\_D incurs a much higher deprivation cost. When more weight is assigned to the deprivation cost then the total cost for RH\_2SSP is significantly lower: 25.49\% lower when the weight is set to 5 and 43.66\% lower when the weight is set to 10. This is due to the inflexibility of RH\_D to handle fluctuating demand and supply. In other words, RH\_D ends up with having to use the more expensive air transportation to satisfy demand. RH\_2SSP provides more robust solutions to handle the uncertainty in post-hurricane logistics operations compared to the RH\_D, reflecting the value of stochastic programming in the rolling-horizon setting.


\subsubsection{Value of Rolling Horizon Approach}
\label{sec:4.5.2}

To investigate the value of the rolling horizon approach we compared the out-of-sample performances of the S\_2SSP and RH\_2SSP approaches. Table \ref{tab:basic_results} summarizes the results where the numbers in the parentheses show the percentages of the corresponding costs in the total cost. For example, for RH\_2SSP, the logistics cost was about 89.1\% of the total cost and the deprivation cost was about 10.9\% of the total. Note that the air transportation cost reported in row four is already included in the logistics cost, but it is presented separately to highlight how differently the two approaches utilize this more expensive transportation mode. As can be seen from Table \ref{tab:basic_results}, on average, the RH\_2SSP resulted in a total cost that is about 11.56\% lower compared to the S\_2SSP. Furthermore, RH\_2SSP resulted in a much lower deprivation cost, whereas the logistics cost is higher. The key reason for a higher deprivation cost incurred by S\_2SSP is that, under each sample path in the out-of-sample test, the deprivation cost is determined by state variables $U$, the opening cost of the SAs is determined by state variables $x$, and the inventory handling cost is determined by state variables $V$, all of which are fixed. The S\_2SSP can only aim to minimize the transportation cost in each stage of the planning horizon by optimizing the local control variables. On the other hand, the RH\_2SSP is rebuilt and re-solved at each roll of each sample path. Thus, RH\_2SSP better exploits the ``local information'' on that sample path. Additionally, we observe that the air transportation cost is much less in the RH\_2SSP, but relatively higher in the S\_2SSP. This is because during the out-of-sample test, the ground transportation is insufficient to address the fluctuating demand (and supply limits) while keeping the values for the state variables (such as the inventory levels) set a priori by the S\_2SSP. In addition, we observe that the RH\_2SSP is computationally much faster: the average solution time was 85 seconds compared to 2,716 seconds for the S\_2SSP. This is because we consider only 10 scenarios for each roll of RH\_2SSP, whereas the S\_2SSP requires solving a 2SSP model with 100 scenarios. The relatively shorter computational time and higher solution quality given by the RH\_2SSP can be attributed to the fact that the RH\_2SSP allows a relatively small number of scenarios in each roll to provide a sufficiently good policy (see our stability test results in Section~\ref{sec:4.5.2}), while in the static case, as we discuss in Section~\ref{sec:4.5.1}, we need to incorporate a relatively large number of scenarios to obtain solutions of acceptable quality. 

\begingroup
\setlength{\tabcolsep}{1pt} 
\begin{table}[htbp]
\footnotesize
\caption{Summary results for RH\_2SSP and S\_2SSP Approaches}
\label{tab:basic_results}
\centering
\begin{tabular}{>{\raggedright}p{0.25\textwidth}>{\raggedleft}p{0.1\textwidth}>{\raggedright}p{0.1\textwidth}>{\raggedleft}p{0.1\textwidth}>{\raggedright}p{0.1\textwidth}>{\raggedleft\arraybackslash}p{0.27\textwidth}}
\hline
 & \multicolumn{2}{c}{RH\_2SSP} & \multicolumn{2}{c}{S\_2SSP} & RH\_2SSP vs. S\_2SSP \\
\hline
\_Total (\$) & 195,564 &  & 221,135 &  & -11.56\%\\
\_Log (\$) & 174,299 & (89.1\%) & 127,493 & (57.6\%) & 36.71\%\\
\_Dep (\$) & 21,265 & (10.9\%) & 93,642 & (42.4\%) & -77.30\%\\
\_Air (\$) & 1 & ($\sim 0\%$) & 2,615 & ($\sim 0\%$) & -99.96\%\\
\_compT (seconds) & 85 & \multicolumn{1}{c}{-} & 2,716 & \multicolumn{1}{c}{-} & -96.87\%\\
\hline
\end{tabular}
\end{table}
\endgroup

\subsection{Sensitivity Analysis}
\label{sec:4.6}
In this section, we present some sensitivity analysis results for the S\_2SSP and the RH\_2SSP under different baseline input parameters values, deprivation cost weights, and lengths of the planning horizon. 

\subsubsection{Baseline Demand and Supply Limits}
\label{sec:4.6.1}
We first compare results obtained by S\_2SSP and RH\_2SSP under different values for the baseline demand and supply limits. In particular, we multiply these parameters by factors ranging from $0.5$ to $2$, as shown in Table~\ref{tab:base_demand}. We can clearly observe from the table that as the baseline demand and supply limits increase the total cost increases for both approaches. Furthermore, the relative performance of the S\_2SSP compared to RH\_2SSP gets worse: the total cost for the S\_2SSP is about 8.54\% worse than that of the RH\_2SSP when baseline values are set to be the lowest ($0.5$), and this gap increases to 27.38\% as baseline values are set to be the highest ($2$). This suggests that the advantage of the RH\_2SSP approach compared to S\_2SSP is more pronounced in the high demand (and supply limit) situations. 

We also observe that the percentage of the deprivation cost for both approaches increases as the multiplier increases. In fact, this percentage suddenly increased to 49.92\% and 84.01\%, respectively, for RH\_2SSP and S\_2SSP when weight was 2. One reason for this increase is that the capacities of the SAs and the PODs are kept the same in these experiments. When the demand becomes high, there is not enough capacity to store the relief items, resulting in high deprivation costs. This observation reflects the importance of an accurate estimation of the post-hurricane demand. Meanwhile, the logistics costs for both S\_2SSP and RH\_2SSP rise as the multiplier increases, but decline when the multiplier is 2. This observation implies that the models tend to deliver smaller amounts of relief items to the affected areas if the regular ground transportation cannot not satisfy the high demand from the PODs and it is too expensive to deliver a large number of relief items via air. One way to avoid this situation is by assigning a high weight to the deprivation cost as we can see in Section~\ref{sec:4.6.2}.

\begin{table}[htbp]
\caption{Performances of S\_2SSP and RH\_2SSP under Different Baseline Demand and Supply Limits}
\label{tab:var_demand_fixcap}       
\centering
\setlength{\tabcolsep}{5pt}
\begin{tabular}{lrrrrrr}
\hline
&\multicolumn{6}{c}{Multiplier} \\ \cline{2-7}
&\multicolumn{1}{c}{0.5}&\multicolumn{1}{c}{0.8} &\multicolumn{1}{c}{1} &\multicolumn{1}{c}{1.2} & \multicolumn{1}{c}{1.5} & \multicolumn{1}{c}{2}\\
\hline\noalign{\smallskip}
RH\_2SSP\_Total (\$) &105,028 &160,707 &195,564 &231,327 & 284,176 & 393,434 \\
RH\_2SSP\_Log (\$) &95,765 &145,308 &174,299 &201,639  & 218,959 &197,019  \\
RH\_2SSP\_Dep (\$)&9,263 &15,399 &21,265 &29,678 & 65,217 &196,414 \\
RH\_2SSP\_Air (\$) &0 &0 &1 &198  &4,993 &21,890\\
RH\_2SSP\_Dep/Total & 8.82\% & 9.58\%& 10.87\% & 12.83\%& 22.95\% & 49.92\%\\
RH\_2SSP\_compT (seconds) &52 &56 &85 &80 & 101 &15\\
\hline
S\_2SSP\_Total (\$) &114,839 &179,789 &221,135 &265,893 &338,571 &541,784   \\
S\_2SSP\_Log (\$)   &85,630  &111,852 &127,493 &144,002 &182,927 &86,640   \\
S\_2SSP\_Dep (\$) &29,209 &67,937  &93,642 &121,891 &155,644 &455,144  \\
S\_2SSP\_2SSP\_Air (\$) &1,467 &2,066 &2,615 &1,960 &1,021 & 4,728 \\
S\_2SSP\_Dep/Total & 25.43\% & 37.79\%& 42.26\% & 45.84\%& 45.97\% & 84.01\%\\
S\_2SSP\_compT (seconds) &1,608 &2,956 &2,716 &3,975 & 4,017 &314 \\
\hline
Improv\_P &8.54\% & 10.61\%& 11.56\%& 13.00\% & 16.07\%  &27.38\%\\
\noalign{\smallskip}\hline
\end{tabular}
\end{table}
\subsubsection{Relative Weight of the Deprivation Cost in the Objective}
\label{sec:4.6.2}
We next present the performances of RH\_2SSP and S\_2SSP under various relative weights of the deprivation cost in the objective function. Specifically, we multiply the $\lambda_{k}$ parameter, which is used as the  weight of the deprivation cost in the objective, with the following factors: 0.25, 0.5, 1, 2 and 5. Intuitively, as this multiplier gets larger, the decision maker puts more emphasis on minimizing the deprivation cost instead of minimizing the transportation cost. Table~\ref{tab:var_dep_cost} summarizes these results.

We observe from the table that when the weight assigned to the deprivation cost is low (0.25 and 0.5), the deprivation costs obtained by the S\_2SSP approach are equal to the total costs. This means that the solutions generated by S\_2SSP do not ship any relief item simply because the unit deprivation cost is less than the unit transportation cost. On the other hand, in these cases, RH\_2SSP ships relief items and manages to reduce the deprivation cost because of its more flexible decisions thanks to the rolling horizon approach. Specifically, the rolling horizon approach compared to the static approach results in solutions that are 26.59\% and 50.4\% lower in total cost for weights 0.25 and 0.5, respectively. When the weight increases to 2 or 5, the total cost for RH\_2SSP increases slowly while the total cost for S\_2SSP grows rapidly. Thus, we can conclude that S\_2SSP is less effective when the weight assigned to the deprivation cost is too low or too high. 

\begin{table}[htbp]
\caption{Performance under Different Weights of Deprivation Cost}
\label{tab:var_dep_cost}       
\centering
\begin{tabular}{lrrrrr}
\hline
&\multicolumn{5}{c}{Weight} \\ \cline{2-6}
&\multicolumn{1}{c}{0.25} &\multicolumn{1}{c}{0.5} &\multicolumn{1}{c}{1} & \multicolumn{1}{c}{2} & \multicolumn{1}{c}{5}\\
\hline\noalign{\smallskip}
RH\_2SSP\_Total (\$) &51,853 &70,065 &195,564 &209,905 & 223,246   \\
RH\_2SSP\_Log (\$) &39,975 &45,095 &174,299 &199,128  & 211,955  \\
RH\_2SSP\_Dep (\$)&11,878 &24,970 &21,265 &10,776 & 11,291 \\
RH\_2SSP\_Air (\$) &11 &76 &1 &19  &235 \\
RH\_2SSP\_Dep/Total &22.91\% & 35.64\%& 10.87\%& 5.13\%& 5.06\% \\
RH\_2SSP\_compT (seconds) &5 &79 &85 &465 & 89 \\
\hline
S\_2SSP\_Total (\$) &70,636 &141,274 &221,135 &310,659 &492,712   \\
S\_2SSP\_Log (\$)   &0        &0 &127,493 &152,867 &276,867   \\
S\_2SSP\_Dep (\$) &70,636 &141,274  &93,642 &157,791 &215,845   \\
S\_2SSP\_Air (\$) &0 &0 &2,615 &14,496 &103,962   \\
S\_2SSP\_Dep/Total &100\% & 100\%& 42.26\%& 50.79\%& 43.81\% \\
S\_2SSP\_compT (seconds) &139 &3,610 &2,716 &4,475 & 1,257  \\
\hline
Improv\_P &26.59\% & 50.40\%& 11.56\%& 32.44\% & 54.71\% \\
\noalign{\smallskip}\hline
\end{tabular}\\
\end{table}

\subsubsection{Length of the Planning Horizon}
\label{sec:4.6.3}
The results presented so far correspond to a planning horizon of $T=5$ periods. According to \cite{vanajakumari_2016}, the response phase typically continues for 3\---7 days. Nevertheless, we increased $T$ to 10 and 15 to see the impact of the length of the planning horizon on solution quality and computational time. Unfortunately, when $T=15$ the optimization solver takes an excessively long time to solve the 2SSP models in both S\_2SSP and RH\_2SSP. Thus, we imposed the following stopping criteria and reported the best (not necessarily optimal) solution found at termination. When $T$ is 5 or 10, we terminated the process after the optimality gap was less than 0.01\% or the computational time reached 24 hours (whichever occurred first) for both S\_2SSP and RH\_2SSP. When $T$ is 15, we set the optimality gap threshold to be 1\% for both approaches, but we set the computational time limit to 24 hours for S\_2SSP and 15 minutes for each roll of RH\_2SSP.

\begin{figure}[htbp]
\begin{center}
  \includegraphics[scale=0.3]{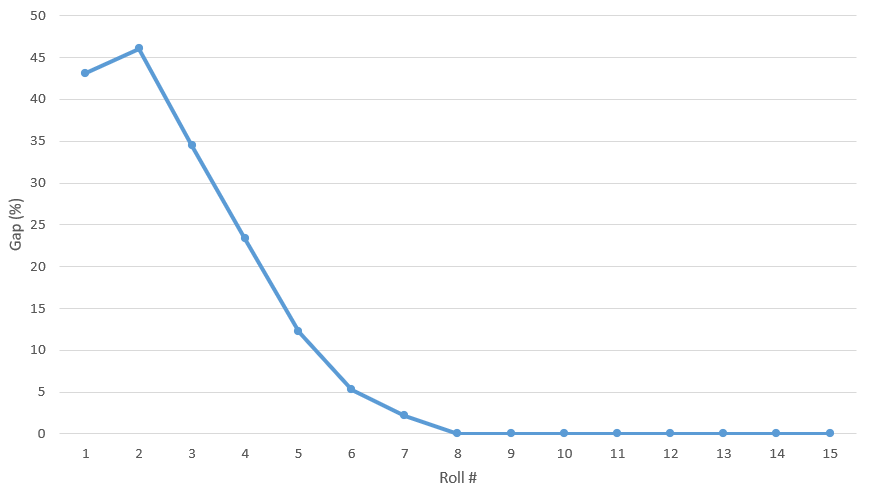}
\vspace{-0cm}
\caption{Illustration of Optimality Gaps of Each Roll of RH\_2SSP When $T = 15$}
\label{fig:6}       
\end{center}
\end{figure}

The results are summarized in Table~\ref{tab:length_th}. The ``NA" in the table indicates that the optimization solver could not yield the desired optimality gap within the specified time limit. For example, for S\_2SSP, when $T=5$ the algorithm terminated by satisfying the optimality gap criterion after 2,716 seconds, whereas when $T=10$($T=15$) the algorithm terminated after reaching the time limit of $24$ hours with an optimality gap of 4.04\%(39.67\%). On the other hand, for RH\_2SSP, when $T=5$($T=15$) the algorithm terminated after 85(3,709) seconds and the optimality gap at each roll was within the threshold of 0.01\%. When $T=15$, the optimality gap threshold for RH\_2SSP was not reached within the time limit in each roll. As can be seen in Figure \ref{fig:6}, each of the first seven rolls took 15 minutes and the optimality gap was over 40\% initially and dropped to about 2\% by the seventh roll. However, the optimality gap threshold was reached for each of the remaining rolls in RH\_2SSP. 

The performances of RH\_2SSP against S\_2SSP on different length of the planning horizon once again demonstrate the flexibility of the rolling horizon approach. Even if bad decisions are made in the first few rolls, RH\_2SSP still ends up with high quality solutions because a new 2SSP model is solved at each roll based on updated information. On the other hand, the quality of S\_2SSP solutions appears to deteriorate as $T$ increases, partly because of the inflexibility in the corresponding decisions as well as the suboptimality in the solutions as a result of the computational time limits.

\begin{table}[htbp]
\caption{Results for Different Lengths of the Planning Horizon }
\label{tab:length_th}       
\centering
\begin{tabular}{lrrr}
\hline
&\multicolumn{3}{c}{Planning Horizon ($T$)} \\ \cline{2-4}
&\multicolumn{1}{c}{5} &\multicolumn{1}{c}{10} &\multicolumn{1}{c}{15} \\
\hline\noalign{\smallskip}
RH\_2SSP\_Total (\$) &195,564 & 386,914 & 584,564 \\
RH\_2SSP\_Log (\$) &174,299 & 331,134& 530,053  \\
RH\_2SSP\_Dep (\$)&21,265 & 55,780& 54,511 \\
RH\_2SSP\_compT (sec) &85 & 3,709& NA\\
\hline
S\_2SSP\_Total (\$) &221,135 & 440,266& 661,084  \\
S\_2SSP\_Log (\$) &127,493 &  218,119& 306,916 \\
S\_2SSP\_Dep (\$) &93,642 &  222,147& 354,168\\
S\_2SSP\_compT (sec) &2,716 & NA&  NA \\
S\_2SSP\_Gap & 0.01\%& 4.04\%&39.67\%\\
\noalign{\smallskip}\hline
\end{tabular}
\end{table}

\section{Conclusions}
\label{sec:5}
We developed static and rolling horizon two-stage stochastic programming approaches for the post-hurricane humanitarian relief logistics planning problem. Our models incorporate both the logistics cost and the social cost via a new definition of the deprivation cost to handle the fluctuating demand. More specifically, we developed a MIP model for the problem and proposed several approaches for implementing this model: deterministic vs. stochastic, and static vs. rolling horizon. We combined a stochastic model that describes the evolution of the post-hurricane logistics system status and the MIP model to formulate a static two-stage stochastic programming model. Furthermore, we integrated the two-stage stochastic programming model within a rolling horizon framework to capture the evolving uncertainty over time. Finally, we conducted extensive computational experiments based on a logistics network of the City of New Orleans. After analyzing our results, we demonstrated the value of stochastic programming solutions compared to their deterministic counterparts in dealing with evolving system state uncertainty arising in post-hurricane relief logistics operations. We also observed that the rolling horizon framework is able to achieve a more cost-effective logistics operation plan than the one provided by a static approach. The advantages of the rolling horizon approach can be seen in terms of both the solution quality and the computational time. 

We identified several directions for future study. One such direction would be to extend the proposed stochastic lookahead framework to address the more challenging vehicle routing and crew scheduling problems for post-hurricane humanitarian relief logistics planning problems. Another possible direction is to validate the presented stochastic lookahead framework using more realistic historical data and real-time data on weather forecasts and hazard analysis. Finally, it would be interesting to compare the proposed rolling horizon two-stage stochastic programming model with a multistage stochastic programming model and analyze their performances in humanitarian relief logistics operation planning.


%
%


\bibliographystyle{spbasic}      

%
%

\newpage
\textbf{Appendix: The static two-stage stochastic programming formulation for the post-hurricane relief logistics planning problem}\\
\begin{subequations}\label{all}
\begin{align}
    Min \quad z = &\sum\limits_{t\in \mathcal{T}} \Big( \sum\limits_{i\in \mathcal{L}}\eta_{i}y_{it} + \sum\limits_{i\in \mathcal{L}}\sum\limits_{k\in \mathcal{K}} \zeta_{k}V_{ikt} \Big) \nonumber\\
    &+ \sum\limits_{\omega \in \Omega} P^{\omega} \sum\limits_{k\in \mathcal{K}}\sum\limits_{t\in \mathcal{T}} \Big(\sum\limits_{i\in L} B^{g}_{i} g^{\omega}_{ikt} + \sum\limits_{i\in \mathcal{L}\cup \mathcal{S}} B^{h}_{i} h^{\omega}_{ikt} + \sum\limits_{i\in \mathcal{L}\cup \mathcal{S}}\sum\limits_{j\in \mathcal{L}\cup \mathcal{S}, j\neq i} B_{ij} f^{\omega}_{ijkt}\Big) \nonumber\\
    &+ \sum\limits_{\omega \in \Omega} P^{\omega} \Big(\sum\limits_{i\in \mathcal{S}}\sum\limits_{k\in \mathcal{K}}\sum\limits_{t\in \mathcal{T}\setminus\{0\}}dep[w, i, k, U_{ik(t-1)}, t-1]\alpha_{ikt} \nonumber\\
    &+ \sum\limits_{i\in \mathcal{S}}\sum\limits_{k\in \mathcal{K}}dep[w, i, k, U_{ikT}, T]\Big) \label{all-obj}\\
    s.t.\quad V_{ik0}, U_{ik0}= &0, \label{all_ini} \qquad\qquad\forall i\in \mathcal{S},k \in \mathcal{K}\\
    V_{ik0}= & g^{\omega}_{ik0}, \label{all_V_init} \qquad\qquad \forall\: i\in \mathcal{L}, k \in \mathcal{K}, \omega \in \Omega\\
    h^{\omega}_{ik0} = &0, \label{all_h_init} \qquad\qquad \forall\: i\in \mathcal{L}\cup S, k \in \mathcal{K}, \omega \in \Omega\\
    \sum \limits_{{\scriptscriptstyle j\in \mathcal{L} \cup \mathcal{S}, j \neq i }}f^{\omega}_{ijk0}= & 0,  \label{all_f_init} \qquad\qquad \forall\: i \in \mathcal{L} \cup \mathcal{S}, k \in \mathcal{K}, \omega \in \Omega\\
    x_{it}= & \sum^{t}_{t^{'} = 0}y_{it^{'}},  \label{all_x_def} \qquad\qquad \forall i\in \mathcal{L}, t\in \mathcal{T}\\
    \sum_{t\in \mathcal{T}}y_{it} \leq & 1, \label{all_y_def} \qquad\qquad \forall i\in \mathcal{L}\\
    V_{ikt} \leq & \phi_{ik} x_{it}, \label{all_SA_cap} \qquad\qquad \forall i\in \mathcal{L}, k \in \mathcal{K}, t\in \mathcal{T}\\
    g^{\omega}_{ikt} + h^{\omega}_{ikt}\leq & V_{ikt}, \quad \label{all_supply_cap} \qquad\qquad \forall\: i\in \mathcal{L}, k \in \mathcal{K}, t\in \mathcal{T}, \omega \in \Omega\\
    \sum\limits_{i\in \mathcal{L}}g^{\omega}_{ikt}\leq & R^{\omega}_{kt}, \quad \label{all_gsupply_cap} \qquad\qquad \forall\: k \in \mathcal{K}, t\in \mathcal{T}, \omega \in \Omega\\
    V_{ikt}+ \sum\limits_{{\scriptscriptstyle j\in \mathcal{L} \cup \mathcal{S}, j \neq i}}f^{\omega}_{ijkt} = & V_{ik(t-1)}+ g^{\omega}_{ikt} + h^{\omega}_{ikt} \nonumber\\
    &+ \sum\limits_{{\scriptscriptstyle j\in \mathcal{L} \cup \mathcal{S}, j \neq i}}f^{\omega}_{jikt}, \label{all_SA_flow} \qquad\qquad \forall\: i \in \mathcal{L}, k \in \mathcal{K}, t\in \mathcal{T}\setminus\{0\}, \omega \in \Omega\\
    U_{ikt}= & (1 - \alpha_{ikt})(U_{ik(t-1)}+1),  \label{all_U_def} \qquad\qquad \forall i \in \mathcal{S}, k \in \mathcal{K}, t\in \mathcal{T}\setminus\{0\}\\
    V_{ikt}\leq & \alpha_{ikt} \phi_{ik}, \label{all_POD_cap_satis} \qquad\qquad \forall i \in \mathcal{S}, k \in \mathcal{K}, t\in \mathcal{T}\setminus\{0\}\\
    \sum\limits_{{\scriptscriptstyle j\in \mathcal{L} \cup \mathcal{S}, j \neq i }}f^{\omega}_{ijkt}\leq & V_{ikt}, \label{all_f_satis1} \qquad\qquad \forall\: i \in \mathcal{S}, k \in \mathcal{K}, t\in \mathcal{T}\setminus\{0\}, \omega \in \Omega \\
    \sum\limits_{{\scriptscriptstyle j\in \mathcal{L} \cup \mathcal{S}, j \neq i }}f^{\omega}_{jikt}\leq & V_{ikt}, \label{all_f_satis2} \qquad\qquad \forall\: i \in \mathcal{S}, k \in \mathcal{K}, t\in \mathcal{T}\setminus\{0\}, \omega \in \Omega\\
    h^{\omega}_{ikt}\leq & V_{ikt}, \label{all_h_satis} \qquad\qquad \forall\: i \in \mathcal{S}, k \in \mathcal{K}, t\in \mathcal{T}\setminus\{0\}, \omega \in \Omega\\
    V_{ikt} + \sum\limits_{{\scriptscriptstyle j\in \mathcal{L} \cup \mathcal{S} \cup \{n\}, j \neq i}}f^{\omega}_{ijkt}=&V_{ik(t-1)}- \alpha_{ikt}  D^{\omega}_{ikt} + h^{\omega}_{ikt} \nonumber\\
    &+ \sum\limits_{{\scriptscriptstyle j\in \mathcal{L} \cup \mathcal{S}, j \neq i}}f^{\omega}_{jikt}, \label{all_POD_flow} \qquad\qquad \forall\: i \in \mathcal{S}, k \in \mathcal{K}, t\in \mathcal{T}\setminus\{0\}, \omega \in \Omega\\
    x_{it},y_{it}\in &\{0,1\}, \label{all_xy_domain} \qquad\qquad \forall i\in \mathcal{L}, t\in \mathcal{T}\\
    \alpha_{ikt}\in &\{0,1\}, \label{all_alpha_domain} \qquad\qquad \forall i \in \mathcal{S}, k \in \mathcal{K}, t\in \mathcal{T}\setminus\{0\}\\
    U_{ikt}\in & \mathbb{Z^+}, \label{all_U_domain} \qquad\qquad \forall i \in \mathcal{S}, k \in \mathcal{K}, t\in \mathcal{T}\\
    V_{ikt}\geq &0, \label{all_V_domain} \qquad\qquad \forall i \in \mathcal{L} \cup S, k \in \mathcal{K}, t\in \mathcal{T}\\
    g^{\omega}_{ikt}, h^{\omega}_{ikt}\geq &0,  \label{all_g&h_domain} \qquad\qquad \forall\: i \in \mathcal{L}, k \in \mathcal{K}, t\in \mathcal{T}, \omega \in \Omega\\
    f^{\omega}_{ijkt}\geq &0, \label{all_f_domain} \qquad\qquad  \forall\: i, j\in \mathcal{L}\cup \mathcal{S} \cup \{n\}, k \in \mathcal{K}, t\in \mathcal{T}, \omega \in \Omega
\end{align}
\end{subequations}

\end{document}